\magnification=1200
\input amstex
\documentstyle{amsppt}

\def\wt#1{\widetilde {#1}}
\def\ssm {\smallsetminus}
\def\SO {\Cal O}
\def\SD {\Cal D}
\def\SE {\Cal E}
\def\SL {\Cal L}
\def\SF {\Cal F}
\def\SH {\Cal H}
\def\SP {\Cal P}
\def\SQ {\Cal Q}
\def\SR {\Cal R}
\def\SU {\Cal U}
\def\SW {\Cal W}

\def\bwq {\bold{\widetilde Q}}

\topmatter
\title Factorization of generalized Theta functions at reducible case  
\endtitle
\author Xiaotao Sun   \endauthor
\address Institute of Mathematics, Academy of Mathematics and 
Systems Sciences, Chinese Academy of Sciences 
Beijing 100080, China \endaddress
\email xsun$\@$math08.math.ac.cn\endemail
\thanks This research was done when the author visited Universit{\"a}t
Essen supported by DFG Forschergruppe "Arithmetik und Geometrie". I thank
all colleagues in Essen for their hospitality.  \endthanks
\endtopmatter
\document

\heading   Introduction  \endheading

One of the problems in algebraic geometry motivated by conformal field theory
is to study the behaviour of moduli space of semistable parabolic bundles on curve and its generalized theta functions when the curve degenerates to a singular curve. Let $X$ be a smooth projective curve of genus $g$, and $\SU_X$
be the moduli space of semistable parabolic bundles on $X$, one can define 
canonically an ample line bundle $\Theta_{\SU_X}$ (theta line bundle) on $\SU_X$ and the global sections $H^0(\Theta^k_{\SU_X})$ are called generalized theta
functions of order $k$. These definitions can be extended to the case of
singular curve. Thus, when $X$ degenerates to a singular curve $X_0$, one may
ask the question how to determine $H^0(\Theta^k_{\SU_{X_0}})$ by generalized theta functions associated with the normalization $\wt X_0$ of $X_0$.
The so called fusion rules suggest that when $X_0$ is a nodal curve the space
$H^0(\Theta^k_{\SU_{X_0}})$ decomposes into a direct sum of spaces of generalized
theta functions on moduli spaces of bundles over $\wt X_0$ with new
parabolic structures at the preimages of nodes. These factorizations and Verlinde 
formula were treated by many mathematicans from various points of view. It is
obviously beyond my ability to give a complete list of contributions. According
to [Be], there are roughly two approachs: infinite and finite. I understand that
those using stacks and loop groups are infinite approach, and working in the
category of schemes of finite type is finite approach. Our approach here should
be a finite one.

When $X_0$ is irreducible with one node, a factorization theorem was proved in
[NR] for rank two and generalized to arbitrary rank in [Su]. By this factorization, one can principally reduce the computation of generalized theta
functions to the case of genus zero with many parabolic points. In order to have
an induction machinery for the number of parabolic points, one should also prove
a factorization when $X_0$ has rwo smooth irreducible components intersecting
at a node $x_0$. This was done for rank two in [DW1] and [DW2] by analytic
method. In this paper, we adopt the approach of [NR] and [Su] to prove a
factorization theorem for arbitrary rank in the reducible case.   

Let $I=I_1\cup I_2\subset X$ be a finite set of points and $\SU_X^I$ the
moduli space of semistable parabolic bundles with parabolic structures at
points $\{x\}_{x\in I}$. When $X$ degenerates to $X_0=X_1\cup X_2$ and
points in $I_j$ ($j=1,2$) degenerate to $|I_j|$ points 
$x\in I_j\subset X_j\ssm\{x_0\}$, 
we have to construct a degeneration $\SU_{X_0}:=\SU^{I_1\cup I_2}_{X_1\cup X_2}$
of $\SU^I_X$ and theta line bundle $\Theta_{\SU_{X_0}}$ on it.
Fix a suitable ample line bundle $\SO(1)$ on $X_0$, we construct the degeneration as a moduli space of `semistable' parabolic torsion free sheaves
on $X_0$ with parabolic structures at points $x\in I_1\cup I_2$, and define 
the theta line bundle $\Theta_{\SU_{X_0}}$ on it. Our main observation here
is that we need a `new semistability' (see definition 1.3) to construct the
correct degeneration of $\SU^I_X.$ But in whole paper, this `new semistability' is simply called semistable. It should not cause any confusion since our
`new semistability' coincides with Seshadri's semistability in [Se] when $I=\emptyset$, and coincides with the semistability of [NR] when $X_0$ is 
irreducible. 

Let $\pi:\wt X_0\to X_0$ be the normalization of $X_0$ and
$\pi^{-1}(x_0)=\{x_1,x_2\}$. Then for any $\mu=(\mu_1,\cdots,\mu_r)$  with $0\le\mu_r\le\cdots\le\mu_1\le k-1,$ we can define $\vec a(x_j)$, $\vec n(x_j)$ and
$\alpha_{x_j}$ ($j=1,2$) by using $\mu$ (see Notation 3.1). Let 
$$\SU^{\mu}_{X_j}:=\SU_{X_j}(r,\chi_j^{\mu}, I_j\cup\{x_j\},
\{\vec n(x),\vec a(x)\}_{x\in I_j\cup\{x_j\}},k)$$
be the moduli space of $s$-equivalence classes of semistable
parabolic bundles $E$ of rank $r$ on $X_j$ and $\chi(E)=\chi_j^{\mu}$,
together with parabolic structures of type 
$\{\vec n(x)\}_{x\in I\cup\{x_j\}}$ and weights 
$\{\vec a(x)\}_{x\in I\cup\{x_j\}}$ at points
$\{x\}_{x\in I\cup\{x_j\}}$, where $\chi_j^{\mu}$ is also defined in Notation 3.1, which may not be integers. Thus we define $\SU^{\mu}_{X_j}$ to be empty 
if $\chi_j^{\mu}$ is not an integer. Let
$$\Theta_{\SU^{\mu}_{X_j}}:=\Theta(k,\ell_j,\{\vec n(x),\vec a(x),\alpha_x\}_{x\in I_j\cup\{x_j\}}, I_j\cup\{x_j\})$$
be the theta line bundle. Then our main result is
\proclaim{Factorization Theorem } There exists a (noncanonical) isomorphism
 $$H^0(\SU_{X_0},\Theta_{\SU_{X_0}})\cong\bigoplus_{\mu}
H^0(\SU^{\mu}_{X_1},\Theta_{\SU^{\mu}_{X_1}})\otimes
H^0(\SU^{\mu}_{X_2},\Theta_{\SU^{\mu}_{X_2}})$$
where $\mu=(\mu_1,\cdots,\mu_r)$ runs through the integers $0\le\mu_r\le\cdots\le
\mu_1\le k-1.$\endproclaim

\S 1 is devoted to construct the moduli space $\SU_{X_0}$ by generalizing 
Simpson's construction, and construct the theta line bundle on it. Then we 
determine the number of irreducible components of the moduli space and proving 
the nonemptyness of them (see Proposition 1.4). In \S 2, we sketch the 
construction of moduli space $\SP$ of generalized parabolic sheaves 
(abbreviated to GPS) and construct an ample line bundle on it. Then we introduce
and study the $s$-equivalence of GPS (see Proposition 2.5), which will be needed
in studying the normalization of $\SU_{X_0}$. In \S 3, we construct and study
the normalization $\SP\to \SU_{X_0}$, then prove the factorization theorem
(Theorem 3.1). As a byproduct, we recover the main results of [NS] (see Corollary 3.1 and Remark 3.1). They have used triples in [NS] instead of GPS.

\remark{Acknowledgements: This work was done during my stay at FB 6 Mathematik
of Universit{\"a}t Essen. I would like to express my hearty thanks
to Prof. H. Esnault and Prof. E. Viehweg for their hospitality and 
encouragements. I was benefited from the stimulating mathematical atmosphere 
they created in their school. Prof. M.S. Narasimhan encouraged me to prove the factorization
theorem at reducible case. I thank him very much for consistent supports.}
\endremark

\heading \S 1 Moduli space of parabolic sheaves \endheading

Let $X_0$ be an reduced projective curve over $\Bbb C$ with two smooth
irreducible components $X_1$ and $X_2$ of genus $g_1$ and $g_2$ meeting
at only one point $x_0$, which is the node of $X_0$. We fix a finite
set $I$ of smooth points on $X_0$ and write $I=I_1\cup I_2$, where 
$I_i=\{x\in I\,|x\in X_i\}$ ($i=1,2$).

\proclaim{Definition 1.1} A coherent $\SO_{X_0}$-module $E$ is called
torsion free if it is pure of dimension $1$, namely, for all nonzero
$\SO_{X_0}$-submodules $E_1\subset E$, the dimension of $Supp(E_1)$
is $1$.\endproclaim

A coherent sheaf $E$ is torsion free if and only if $E_x$ has depth $1$
at every $x\in X_0$ as a $\SO_{X_0,x}$-module. Thus $E$ is locally free
over $X_0\ssm\{x_0\}$.

\proclaim{Definition 1.2} We say that a torsion free sheaf $E$ over $X_0$
has a quasi-parabolic structure of type $\vec n(x)=(n_1(x),\cdots,n_{l_x+1}(x))$
at $x\in I$, if we choose a flag of subspaces
$$E|_{\{x\}}=F_0(E)_x\supset F_1(E)_x\supset\cdots\supset F_{l_x}(E)_x\supset
F_{l_x+1}(E)_x=0$$
such that $n_j(x)=dim(F_{j-1}(E)_x/F_j(E)_x)$. If, in addition, a sequence of
integers called the parabolic weights
$$0 < a_1(x)<a_2(x)<\cdots<a_{l_x+1}(x)< k$$
are given, we call that $E$ has a parabolic structure of type $\vec n(x)$ at 
$x$, with weights $\vec a(x):=(a_1(x),\cdots,a_{l_x+1}(x))$. The sheaf $E$
is also simply called a parabolic sheaf, whose parabolic Euler characteristic is
defined as  
$$par\chi(E):=\chi(E)+\frac{1}{k}\sum_{x\in I}\sum^{l_x+1}_{i=1}n_i(x)a_i(x).$$
\endproclaim

We will fix an ample line bundle $\SO(1)$ on $X_0$ such that 
$deg(\SO(1)|_{X_i})=c_i>0$ ($i=1,2$), for simplicity, we assume that
$\SO(1)=\SO_{X_0}(c_1y_1+c_2y_2)$ for two fixed smooth points $y_i\in X_i$.
For any torsion free sheaf $E$, $P(E,n):=\chi(E(n))$ denote its Hilbert 
polynomial, which has degree $1$. We define the rank of $E$ to be
$$rk(E):=\frac{1}{deg(\SO(1))}\cdot \lim \limits_{n\to\infty}\frac{P(E,n)}
{n}.$$
Let $r_i$ denote the rank of the restriction of $E$ to $X_i$ ($i=1,2$), then
$$P(E,n)=(c_1r_1+c_2r_2)n+\chi(E),\quad r(E)=\frac{c_1}{c_1+c_2}r_1+\frac{c_2}{c_1+c_2}r_2.$$

\proclaim{Notation 1.1} We say that $E$ is a torsion free sheaf of rank $r$ on $X_0$ if $r_1=r_2=r$, otherwise it will be said of rank $(r_1,r_2)$.
We will fix in this paper the parabolic datas $\{\vec n(x)\}_{x\in I}$, 
$\{\vec a(x)\}_{x\in I}$ and the integers: $\chi=d+r(1-g)$, $\ell_1+\ell_2$, $k$, 
$$\alpha:=\{0\le\alpha_x< k-a_{l_x+1}(x)+a_1(x)\}_{x\in I}$$
such that
$$ \sum_{x\in
I}\sum^{l_x}_{i=1}d_i(x)r_i(x)+r\sum_{x\in
I}
\alpha_x+r(\ell_1+\ell_2)=k\chi,\tag{$*$}$$
where $d_i(x)=a_{i+1}(x)-a_i(x)$ and $r_i(x)=n_1(x)+\cdots+n_i(x)$. We will
choose $c_1$ and $c_2$ such that $\ell_1=\frac{c_1}{c_1+c_2}(\ell_1+\ell_2)$
and $\ell_2=\frac{c_2}{c_1+c_2}(\ell_1+\ell_2)$ become integers.
\endproclaim

\proclaim{Definition 1.3} With the fixed parabolic datas in Notation 1.1, and
for any torsion free sheaf $F$ of rank $(r_1,r_2)$, let
$$m(F):= \frac{r(F)-r_1}{k}\sum_{x\in I_1}(a_{l_x+1}(x)+\alpha_x)+
\frac{r(F)-r_2}{k}\sum_{x\in I_2}(a_{l_x+1}(x)+\alpha_x).$$
If $F$ has parabolic structures at points $x\in I$, the modified parabolic
Euler characteristic and slop of $F$ are defined as 
$$par\chi_m(F):=par\chi(F)+m(F),\quad par\mu_m(F):=\frac{par\chi_m(F)}{r(F)}.$$
A parabolic sheaf $E$ is called semistable (resp. stable) for
$(k,\alpha,\vec a)$ if, for any subsheaf $F\subset E$
such $E/F$ is torsion free, one has, with the induced parabolic structure, 
$$par\chi_m(F)\le \frac{par\chi_m(E)}{r(E)}r(F)\quad (resp.<).$$
\endproclaim

\remark{Remark 1.1} The above semistability is independent of the choice of
$\alpha$ and coincides with Seshadri's semistability of parabolic torsion free
sheaves when the curves are irreducible.\endremark

We will only consider in this section torsion free sheaves 
of rank $r$ with parabolic structures of type $\{\vec n(x)\}_{x\in I}$
and weights $\{\vec a(x)\}_{x\in I}$ at points $\{x\}_{x\in I}$, and construct
the moduli space of semistable parabolic sheaves. Let $\SW=\SO_{X_0}(-N)$ and $V=\Bbb C^{P(N)}$, we consider the Quot scheme
$$Quot(V\otimes\SW,P)(T)=\left\{
\aligned &\text{$T$-flat quotients $V\otimes\SW\to E\to 0$ over}\\ 
&\text{$X_0\times T$ with Hilbert polynomial $P$}\endaligned\right\},$$
and let $\bold Q\subset Quot(V\otimes\SW,P)$ be the open set
$$\bold Q(T)=\left\{\aligned&\text{$V\otimes\SW\to E\to 0$, with 
$R^1p_{T*}(E(N))=0$ and}\\&\text{$V\otimes\SO_T\to p_{T*}E(N)$ 
induces an isomorphism}\endaligned\right\}.$$ Thus we can assume (Lemma 20 of [Se], page 162) that $N$ is chosen large enough so that every semistable
parabolic torsion free sheaf with Hilbert polynomial $P$ and parabolic structures of type $\{\vec n(x)\}_{x\in I}$, weights $\{\vec a(x)\}_{x\in I}$
at points $\{x\}_{x\in I}$ appears as a quotient corresponding to a point of
$\bold Q$. 

Let $V\otimes\SW\to \SF\to 0$ be the universal quotient over $X_0\times\bwq$
and $\SF_x$ be the restrication of $\SF$ on $\{x\}\times\bwq\cong\bwq$.
Let $Flag_{\vec n(x)}(\SF_x)\to\bwq$ be the relative flag scheme of type
$\vec n(x)$, and 
$$\SR=\underset x\in I \to{\times_{\bwq}}Flag_{\vec n(x)}(\SF_x)\to\bwq$$
be the product over $\bwq$, where $\bwq$ is the closure of $\bold Q$ in the
Quot scheme.  

A (closed) point $(p,\{p_{r(x)},p_{r_1(x)},...,p_{r_{l_x}(x)}\}_{x\in I})$ of $\SR$ by definition
is given by a point $V\otimes\SW@>p>>E\to 0$ of the Quot scheme, together with
quotients $$\{V\otimes\SW@>p_{r(x)}>>Q_{r(x)},
V\otimes\SW@>p_{r_1(x)}>>Q_{r_1(x)},...,
V\otimes\SW@>p_{r_{l_x(x)}}>>Q_{r_{l_x}(x)}\}_{x\in I},$$
where $r_i(x)=dim(E_x/F_i(E)_x)=n_1(x)+\cdots+n_i(x)$, and $Q_{r(x)}:=E_x,$
$Q_{r_1(x)}:=E_x/F_1(E)_x, ..., Q_{r_{l_x}(x)}:=E_x/F_{l_x}(E)_x$.
For large enough $m$, we have a $SL(V)$-equivariant embedding
$$\Cal R\hookrightarrow \bold G=Grass_{P(m)}(V\otimes W_m)\times\bold {Flag},$$
where $W_m=H^0(\SW(m))$, and $\bold{Flag}$ is defined to be
$$\bold{Flag}=\prod_{x\in I}\{Grass_{r(x)}(V)\times Grass_{r_1(x)}(V)
\times\cdots\times
Grass_{r_{l_x}(x)}(V)\},$$
which maps a point $(p,\{p_{r(x)},p_{r_1(x)},...,p_{r_{l_x}(x)}\}_{x\in I})$
of $\SR$ to the point
$$(V\otimes W_m@>g>>U,\{V@>g_{r(x)}>>U_{r(x)},
V@>g_{r_1(x)}>>U_{r_1(x)},...,
V@>g_{r_{l_x}(x)}>>U_{r_{l_x}(x)}\}_{x\in I})$$
of $\bold G$, where $g:=H^0(p(m))$, $U:=H^0(E(m))$, $g_{r(x)}:=H^0(p_{r(x)}(N)),$ $U_{r(x)}:=H^0(Q_{r(x)}),$ and $g_{r_i(x)}:=H^0(p_{r_i(x)}(N))$, $U_{r_i(x)}:=H^0(Q_{r_i(x)})$ ($i=1,...,l_x$).

For any rational number $\ell$ satisfying $c_i\ell=\ell_i+c_ikN$ ($i=1,2$), we give $\bold G$ the polarisation (using
the obvious notation):
$$\frac{\ell}{m-N}\times\prod_{x\in
I}\{\alpha_x, d_1(x),\cdots,
d_{l_x}(x)\},$$
and we have a straightforward generalisation of [NR, Proposition A.6] whose
proof we omit:
\proclaim{Proposition 1.1} A point 
$(g,\{g_{r(x)},g_{r_1(x)},...,g_{r_{l_x}(x)}\}_{x\in I})\in \bold G$ is stable
(respectively, semistable) for the action of $SL(V)$, with respect to the 
above polarisation (we refer to this from now on as GIT-stability), iff for all
nontrivial subspaces $H\subset V$ we have (with $h=dim H$)
$$\aligned\frac{\ell}{m-N}(hP(m)&-P(N)dim\,g(H\otimes W_m))+
\sum_{x\in I}\alpha_x(rh-P(N)dim\,g_{r(x)}(H))\\
&+\sum_{x\in I}\sum^{l_x}_{i=1}d_i(x)(r_i(x)h-P(N)dim\,g_{r_i(x)}(H))
<(\le)\,0.\endaligned$$\endproclaim

\proclaim{Notation 1.2} Given a point  $(p,\{p_{r(x)},p_{r_1(x)},...,p_{r_{l_x}(x)}\}_{x\in I})\in\SR$, and a subsheaf
$F$ of $E$ we denote the image of $F$ in $Q_{r_i(x)}$ (respectively, in $Q_{r(x)}$) by $Q_{r_i(x)}^F$ (respectively, by $Q_{r(x)}^F$). Similarly, given a quotient $E@>T>>\Cal G\to 0$, set $Q_{r_i(x)}^{\Cal G}:=Q_{r_i(x)}/Im(ker(T))$
(respectively, $Q_{r(x)}^{\Cal G}:=Q_{r(x)}/Im(ker(T))$).\endproclaim

\proclaim{Proposition 1.2} Suppose 
$(p,\{p_{r(x)},p_{r_1(x)},...,p_{r_{l_x}(x)}\}_{x\in I})\in\SR$ is a point
such that $E$ is torsion free. Then $E$ is stable (respectively, semistable)
iff for every subsheaf $0\neq F\neq E$ we have
$$\aligned\frac{\ell}{m-N}&(\chi(F(N))P(m)-P(N)\chi(F(m)))+
\sum_{x\in I}\alpha_x(r\chi(F(N))-P(N)h^0(Q_{r(x)}^F))\\
&+\sum_{x\in I}\sum^{l_x}_{i=1}d_i(x)(r_i(x)\chi(F(N))-P(N)h^0(Q_{r_i(x)}^F))
<(\le)\,0.\endaligned$$\endproclaim

\demo{Proof} For any subsheaf $F$ let $LHS(F)$ denote the left-hand side of
above inequality. Assumme first that $E/F$ is torsion free and $F$ is of rank
$(r_1,r_2)$, thus $h^0(Q^F_{r(x)})=r_i$ for $x\in I_i$ ($i=1,2$) and
$\chi(F(m))=(c_1r_1+c_2r_2)(m-N)+\chi(F(N))$, 
$$h^0(Q^F_{r_i(x)})=dim(F_x/F_x\cap F_i(E)_x).$$
Let $n^F_i(x):= dim(F_x\cap F_{i-1}(E)_x/F_x\cap F_i(E)_x)$, and note that
$$\sum_{x\in I}\sum^{l_x}_{i=1}d_i(x)r_i(x)=r\sum_{x\in I}a_{l_x+1}(x)
-\sum_{x\in I}\sum^{l_x+1}_{i=1}a_i(x)n_i(x),$$
$$\aligned\sum_{x\in I}\sum^{l_x}_{i=1}d_i(x)dim(F_x/F_x\cap F_i(E)_x)&=
r_1\sum_{x\in I_1}a_{l_x+1}(x)+r_2\sum_{x\in I_2}a_{l_x+1}(x)\\
&-\sum_{x\in I}\sum^{l_x+1}_{i=1}a_i(x)n^F_i(x),\endaligned$$
we have
$$\aligned LHS(F)&=\left(\sum_{x\in I}\sum^{l_x}_{i=1}d_i(x)r_i(x)+
r\sum_{x\in I}\alpha_x+rc_1\ell+rc_2\ell\right)(\chi(F)-\frac{r(F)}{r}\chi)\\
&+P(N)\left( r(F)\sum_{x\in I}\alpha_x+ \frac{r(F)}{r}\sum_{x\in I}
\sum^{l_x}_{i=1}d_i(x)r_i(x)\right)\\
&-P(N)\left(r_1\sum_{x\in I_1}\alpha_x+r_2\sum_{x\in I_2}\alpha_x
+\sum_{x\in I}\sum^{l_x}_{i=1}d_i(x)dim(F_x/F_x\cap F_i(E)_x)
\right)\\
&=kP(N)\left(par\chi_m(F)-\frac{r(F)}{r}par\chi_m(E)\right).
\endaligned$$
Thus the inequality implies the (semi)stability of $E$, and the (semi)stability
of $E$ implies the inequality for subsheaves $F$ such that $E/F$ torsion free.

Suppose now that $E$ is (semi)stable and $F$ any nontrivial subsheaf, let
$\tau$ be the torsion of $E/F$ and $F'\subset E$ such that $\tau=F'/F$ and
$E/F'$ torsion free. Then we have $LHS(F')\le 0$ and, if we write $\tau=\tilde\tau+\sum_{x\in I}\tau_x,$ then
$$\aligned &LHS(F)-LHS(F')= -(c_1+c_2)r\ell h^0(\tau)\\&-
\sum_{x\in I}\alpha_x\left(rh^0(\tau)+P(N)(h^0(Q^F_{r(x)})-h^0(Q^{F'}_{r(x)}))
\right)\\&-\sum_{x\in I}\sum^{l_x}_{i=1}d_i(x)\left(r_i(x)h^0(\tau)+P(N)(
h^0(Q^F_{r_i(x)})-h^0(Q^{F'}_{r_i(x)}))\right)\\&=-h^0(\tau) 
\left(\sum_{x\in I}\sum^{l_x}_{i=1}d_i(x)r_i(x)+
r\sum_{x\in I}\alpha_x+rc_1\ell+rc_2\ell\right)-\\&P(N)
\left(\sum_{x\in I}-\alpha_xh^0(\tau_x)+\sum_{x\in I}\sum^{l_x}_{i=1}d_i(x)(h^0(Q^F_{r_i(x)})-h^0(Q^{F'}_{r_i(x)}))\right)\\
&\le-kP(N)h^0(\tau)+P(N)\sum_{x\in I}\alpha_xh^0(\tau_x)+P(N)\sum_{x\in I}\sum^{l_x}_{i=1}d_i(x)h^0(\tau_x)\\
&= -kP(N)h^0(\tilde\tau)-P(N)\sum_{x\in I}(k-\alpha_x-a_{l_x+1}(x)+a_1(x))
h^0(\tau_x)\le 0,
\endaligned$$
where we have used $h^0(Q^F_{r(x)})-h^0(Q^{F'}_{r(x)})=-h^0(\tau_x)$, the assumption about $\{\alpha_x\}$ and
$h^0(Q^F_{r_i(x)})-h^0(Q^{F'}_{r_i(x)})\ge -h^0(\tau_x)$.\enddemo

\proclaim{Lemma 1.1} There exists $M_1(N)$ such that for $m\ge M_1(N)$ the
following holds. Suppose $(p,\{p_{r(x)},p_{r_1(x)},...,p_{r_{l_x}(x)}\}_{x\in I})\in\SR$ is a point which is GIT-semistable then for all quotients 
$E@>T>>\Cal G\to 0$ we have
$$h^0(\Cal G(N))\ge\frac{1}{k}\left((c_1+c_2)r(\Cal G)\ell+
\sum_{x\in I}\alpha_xh^0(Q^{\Cal G}_{r(x)})+\sum_{x\in I}\sum^{l_x}_{i=1}
d_i(x)h^0(Q^{\Cal G}_{r_i(x)})\right).$$ In particular, $E$ is torsion free and
$V\to H^0(E(N))$ is an isomorphism. \endproclaim

\demo{Proof} Let $H=ker\{V@>H^0(p(N))>>H^0(E(N))\to H^0(\Cal G(N))\}$, and 
$F\subset E$ the subsheaf generated by $H$. Since all these $F$ are in a bounded
family, $dim\,g(H\otimes W_m)=h^0(F(m))=\chi(F(m))$ for $m$ large enough. Thus
there exists $M_1(N)$ such that for $m\ge M_1(N)$ the inequality of Proposition 1.1 implies (with $h=dim(H)$)
$$\aligned&(c_1+c_2)\ell(rh-r(F)P(N))+\sum_{x\in I}\alpha_x\left(rh-P(N)h^0(Q^F_{r(x)})\right)\\ 
&+\sum_{x\in I}\sum^{l_x}_{i=1}d_i(x)\left(r_i(x)h-P(N)h^0(Q^F_{r_i(x)})\right)
\,\le 0,\endaligned$$
where we used that $g_{r(x)}(H)=h^0(Q^F_{r(x)})$ and $g_{r_i(x)}(H)=h^0(Q^F_{r_i(x)})$. Now using the following inequalities
$$h\ge P(N)-h^0(\Cal G(N)),\quad r-h^0(Q^F_{r(x)})\ge h^0(Q^{\Cal G}_{r(x)}),$$
$$r-r(F)\ge r(\Cal G),\quad r_i(x)-h^0(Q^F_{r_i(x)})\ge 
h^0(Q^{\Cal G}_{r_i(x)}),$$
we get the inequality
$$h^0(\Cal G(N))\ge\frac{1}{k}\left((c_1+c_2)r(\Cal G)\ell+
\sum_{x\in I}\alpha_xh^0(Q^{\Cal G}_{r(x)})+\sum_{x\in I}\sum^{l_x}_{i=1}
d_i(x)h^0(Q^{\Cal G}_{r_i(x)})\right).$$

Now we show that $V\to H^0(E(N))$ is an isomorphism. That it is injective is
easy to see: let $H$ be its kernel, then $g(H\otimes W_m)=0$, $g_{r(x)}(H)=0$
and $g_{r_i(x)}(H)=0$, one sees that $h=0$ from Proposition 1.1. To see it being
surjective, it is enough to show that one can choose $N$ such that $H^1(E(N))=0$
for all such $E$. If $H^1(E(N))$ is nontrivial, then there is a nontrivial 
quotient $E(N)\to L\subset\omega_{X_0}$ by Serre duality, and thus
$$h^0(\omega_{X_0})\ge h^0(L)\ge (c_1+c_2)N+B,$$
where $B$ is a constant independent of $E$, we choose $N$ such that $H^1(E(N))=0$ for all GIT-semistable points.   

Let $\tau=Tor(E)$, $\Cal G=E/\tau$ and applying the above inequality, noting that
$h^0(\Cal G(N))=P(N)-h^0(\tau)$, $h^0(Q^{\Cal G}_{r(x)})=r-h^0(Q^{\tau}_{r(x)})$ and
$h^0(Q^{\Cal G}_{r_i(x)})=r_i(x)-h^0(Q^{\tau}_{r_i(x)}),$ we have
$$kh^0(\tau)\le\sum_{x\in I}(\alpha_x+a_{l_x+1}(x)-a_1(x))h^0(\tau_x),$$
by which one can conclude that $\tau=0$ since $\alpha_x<k-a_{l_x+1}(x)+a_1(x).$
\enddemo

\proclaim{Proposition 1.3} There exist integers $N>0$ and $M(N)>0$ such that
for $m\ge M(N)$ the following is true. A point
$(p,\{p_{r(x)},p_{r_1(x)},...,p_{r_{l_x}(x)}\}_{x\in I})\in\SR$ is GIT-stable
(respectively, GIT-semistable) iff the quotient $E$ is torsion free and a
stable (respectively, semistable) sheaf, the map $V\to H^0(E(N))$ is an isomorphism.\endproclaim

\demo{Proof} If $(p,\{p_{r(x)},p_{r_1(x)},...,p_{r_{l_x}(x)}\}_{x\in I})\in\SR$ is GIT-stable (GIT-semistable), by Lemma 1.1, $E$ is torsion free and
$V\to H^0(E(N))$ is an isomorphism. For any subsheaf $F\subset E$ with
$E/F$ torsion free, let $H\subset V$ be the inverse image of $H^0(F(N))$ and
$h=dim(H)$, we have $\chi(F(N))P(m)-P(N)\chi(F(m))\le hP(m)-P(N)h^0(F(m))$
for $m>N$ (note that $h^1(F(N))\ge h^1(F(m))$). Thus
$$\aligned&\frac{\ell}{m-N}(\chi(F(N))P(m)-P(N)\chi(F(m)))+
\sum_{x\in I}\alpha_x(r\chi(F(N))-P(N)h^0(Q_{r(x)}^F))\\
&+\sum_{x\in I}\sum^{l_x}_{i=1}d_i(x)(r_i(x)\chi(F(N))-P(N)h^0(Q_{r_i(x)}^F))
\le\\
&\frac{\ell}{m-N}(hP(m)-P(N)dim\,g(H\otimes W_m))
+\sum_{x\in I}\alpha_x(rh-P(N)dim\,g_{r(x)}(H))\\
&+\sum_{x\in I}\sum^{l_x}_{i=1}d_i(x)(r_i(x)h-P(N)dim\,g_{r_i(x)}(H))
\endaligned$$
since $g(H\otimes W_m)\le h^0(F(m))$, $g_{r(x)}(H)\le h^0(Q^F_{r(x)})$
and $g_{r_i(x)}(H)\le h^0(Q^F_{r_i(x)})$ (the inequalities are strict when $h=0$). By Proposition 1.1 and Proposition 1.2, $E$ is stable (respectively, semistable) if the 
point is GIT stable (respectively, GIT semistable).

The proof of another direction is similar to [NR], one can prove the similar
Lemma A.9 and Lemma A.12 of [NR] by just modifying notation. 

\enddemo

One can imitate [Se] (Th{\'e}or{\`e}me 12, page 71) to show that
given a semistable parabolic sheaf $E$, there exists a filtration of $E$
$$0=E_{n+1}\subset E_n\subset\cdots\subset E_2\subset E_1\subset E_0=E$$
such that $E_i/E_{i+1}$ ($1\le i\le n$) are stable parabolic sheaves with
the constant slop $par\mu_m(E)$, and the
isomorphic class of semistable parabolic sheaf 
$$gr(E):=\bigoplus_i^nE_i/E_{i+1}$$
is independent of the filtration. Two semistable parabolic sheaves $E$ and $E'$
are called $s$-equivalent if $gr(E)\cong gr(E').$ 

\proclaim{Theorem 1.1} For given datas in Notation 1.1 satisfying $(*)$, there
exists a reduced, seminormal projective scheme 
$$\SU_{X_0}:=\SU_{X_0}(r,\chi,I_1\cup I_2,
\{\vec n(x),\vec a(x),\alpha_x\}_{x\in I},\SO(1), k),$$
which is the coarse moduli space of $s$-equivalence classes of semistable
parabolic sheaves $E$ of rank $r$ and $\chi(E)=\chi$ with parabolic structures of type $\{\vec n(x)\}_{x\in I}$ and weights $\{\vec a(x)\}_{x\in I}$ at points
$\{x\}_{x\in I}$. The moduli space $\SU_{X_0}$ has at most $r+1$ irreducible components.\endproclaim  

\demo{Proof} Let $\SR^{ss}$ ($\SR^s$) be the open set of $\SR$
whose points correspond to semistable (stable) parabolic sheaves on $X_0$.
Then, by Proposition 1.3, the quotient 
$$\varphi: \SR^{ss}\to \SU_{X_0}:=\SR^{ss}//SL(V)$$
exists as a projective scheme. That $\SU_{X_0}$ is reduced and
seminormal follow from the properties of $\SR^{ss}$ (see [Fa] [Se] [Su]).

Consider the dense open set $\SR_0\subset\SR^{ss}$ consists of locally free
sheaves, for each $F\in \SR_0$, let $F_1$ and $F_2$ be the restrictions of
$F$ to $X_1$ and $X_2$, we have
$$0\to F_1(-x_0)\to F\to F_2\to 0.\tag1.1$$
By the semistability of $F$ and $par\chi_m(F_1)+par\chi_m(F_2)=par\chi_m(F)+r$,
we have 
$$\frac{c_1}{c_1+c_2}par\chi_m(F)\le par\chi_m(F_1)\le\frac{c_1}{c_1+c_2}
par\chi_m(F)+r,$$ 
$$\frac{c_2}{c_1+c_2}par\chi_m(F)\le par\chi_m(F_2)\le\frac{c_2}{c_1+c_2}
par\chi_m(F)+r.$$
Let $\chi_1$, $\chi_2$ denote $\chi(F_1)$, $\chi(F_2)$, and $n_j$ denote 
(for $j=1,2$)
$$ \frac{1}{k}\left(\sum_{x\in
I_j}\sum^{l_x}_{i=1}d_i(x)r_i(x)+r\sum_{x\in
I_j}
\alpha_x+r\ell_j\right),\tag1.2$$
we can rewrite the above inequalities into
$$n_1\le\chi_1\le n_1+r,\quad n_2\le\chi_2\le n_2+r.\tag1.3$$
There are at most $r+1$ possible choices of $(\chi_1,\chi_2)$ satisfying
$(1.3)$ and $\chi_1+\chi_2=\chi+r$, each of the choices corresponds an
irreducible component of $\SU_{X_0}$.  
\enddemo 

For any $\chi_1$, $\chi_2$ satisfying (1.3), let $\SU_{X_1}$ (resp. $\SU_{X_2}$)
be the moduli space of semistable parabolic bundles of rank $r$ and Euler
characteristic $\chi_1$ (resp. $\chi_2$), with parabolic structures of type 
$\{\vec n(x)\}_{x\in I_1}$ (resp. $\{\vec n(x)\}_{x\in I_2}$) and weights 
$\{\vec a(x)\}_{x\in I_1}$ (resp. $\{\vec a(x)\}_{x\in I_2}$) 
at points $\{x\}_{x\in I_1}$ (resp. $\{x\}_{x\in I_2}$). Then we have

\proclaim{Proposition 1.4} Suppose that $\SU_{X_1}$ and $\SU_{X_2}$ are not
empty. Then there exists a semistable parabolic vector bundle $E$ on $X_0$,
with parabolic structures of type $\{\vec n(x)\}_{x\in I}$ and weights 
$\{\vec a(x)\}_{x\in I}$ at points $\{x\}_{x\in I}$, such that
$$E|_{X_1}\in \SU_{X_1},\quad E|_{X_2}\in \SU_{X_2}.$$
Moreover, if $n_1\le\chi_1< n_1+r$ and $n_2\le\chi_2< n_2+r$, $E$ is stable
whenever one of $E_1$ and $E_2$ is stable.
\endproclaim

\demo{Proof} For any $E_1\in\SU_{X_1}$ and $E_2\in\SU_{X_2}$, one can glue
them by any isomorphism at $x_0$ into a vector bundle $E$ on $X_0$ with the described
parabolic structures at points $\{x\}_{x\in I}$ such that $E|_{X_1}=E_1$ and $E|_{X_2}=E_2.$ We will show that $E$ is semistable.

For any subsheaf $F\subset E$ of rank $(r_1,r_2)$ such that $E/F$ torsion
free, we have the commutative diagram
$$\CD
0@>>> F_1@>>>F@>>> F_2@>>>  0  \\
@.    @VVV     @VVV        @VVV        @. \\
0@>>> E_1(-x_0)  @>>> E   @>>>  E_2 @>>> 0
\endCD\tag1.4$$ 
where $F_2$ is the image of $F$ under $E\to E_2\to 0$ and $F_1$ is the kernel
of $F\to F_2\to 0.$ One sees easily that $F_1$ and $F_2$ are of rank $(r_1,0)$ and $(0,r_2)$ torsion free sheaves. From the above diagram (1.4),
we have the following equalities 
$$\aligned&\frac{par\chi_m(E)}{r}-\frac{par\chi_m(F)}{r(F)}\\
&=\frac{par\chi_m(E_1(-x_0))}{r}-\frac{par\chi_m(F_1)}{r(F)}+
\frac{par\chi_m(E_2)}{r}-\frac{par\chi_m(F_2)}{r(F)}\\&
= \frac{a_1r_1\cdot par\chi_m(E_1(-x_0))-r\cdot par\chi_m(F_1)+a_2r_2
\cdot par\chi_m(E_1(-x_0))}{r(F)\cdot r}\\&+    
 \frac{a_2r_2\cdot par\chi_m(E_2)-r\cdot par\chi_m(F_2)+a_1r_1
\cdot par\chi_m(E_2)}{r(F)\cdot r}\\
&=\frac{r_1}{r(F)}(par\mu_m(E_1(-x_0))-par\mu_m(F_1))+  \frac{r_2}{r(F)}(par\mu_m(E_2)-par\mu_m(F_2))\\&
+\frac{a_2(r_2-r_1)par\chi_m(E_1(-x_0))+a_1(r_1-r_2)par\chi_m(E_2)}{r(F)\cdot r}\\
&=\frac{r_1}{r(F)}(par\mu(E_1(-x_0))-par\mu(F_1))+  \frac{r_2}{r(F)}(par\mu(E_2)-par\mu(F_2))
\\&+\frac{(r_1-r_2)(\frac{c_1}{c_1+c_2}par\chi_m(E)+r-par\chi_m(E_1))}
{r(F)\cdot r}\endaligned$$
where we used the notation $a_1:=\frac{c_1}{c_1+c_2}$ and 
$a_2:=\frac{c_2}{c_1+c_2},$ the last equality follows that 
$$\frac{m(E_1(-x_0))}{r}-\frac{m(F_1)}{r_1}=0,\quad \frac{m(E_2)}{r}-\frac{m(F_2)}{r_2}=0.$$
Similarly, if we use the following diagram
$$\CD
0@>>> F_2@>>>F@>>> F_1@>>>  0  \\
@.    @VVV     @VVV        @VVV        @. \\
0@>>> E_2(-x_0)  @>>> E   @>>>  E_1 @>>> 0
\endCD$$
we will get the equality
$$\aligned&\frac{par\chi_m(E)}{r}-\frac{par\chi_m(F)}{r(F)}=
\frac{r_2}{r(F)}(par\mu(E_2(-x_0))-par\mu(F_2))\\
&+  \frac{r_1}{r(F)}(par\mu(E_1)-par\mu(F_1))
\\&+\frac{(r_2-r_1)(\frac{c_2}{c_1+c_2}par\chi_m(E)+r-par\chi_m(E_2))}
{r(F)\cdot r}\endaligned$$
Thus we always have the inequality 
$$\frac{par\chi_m(E)}{r}-\frac{par\chi_m(F)}{r(F)}\ge 0$$
and the equality implies that $r_1=r_2$ and $E_1$, $E_2$ are both unstable.
This proves the proposition.
\enddemo

By a family of
parabolic sheaves of rank $r$ and Euler characteristic $\chi$ with parabolic structures of type $\{\vec n(x)\}_{x\in I}$ and weights $\{\vec a(x)\}_{x\in I}$ at points $\{x\}_{x\in I}$ parametrized by $T$, we mean a sheaf $\Cal F$ on
$X_0\times T$, flat over
$T$, and torsion free with rank $r$ and Euler characteristic $\chi$ on
$X_0\times\{t\}$ for every $t\in T$, together with, for each $x\in I$, a
flag 
$$\SF_{\{x\}\times T}=F_0(\SF_{\{x\}\times T})\supset F_1(\SF_{\{x\}\times T}) 
\supset\cdots\supset F_{l_x}(\SF_{\{x\}\times T})\supset 
F_{l_x+1}(\SF_{\{x\}\times T})=0$$
of subbundles of type $\vec n(x)$ and  weights $\vec a(x)$. Let 
$\SQ_{\{x\}\times T,i}$ denote the quotients
$\SF_{\{x\}\times T}/F_i(\SF_{\{x\}\times T})$, then we define a line bundle
$\Theta_{\SF}$ on $T$ to be
$$(detR\pi_T\Cal F)^k\otimes\bigotimes_{x\in I}
\lbrace(det\Cal F_{\{x\}\times
T})^{\alpha_x}\otimes\bigotimes^{l_x}_{i=1}
(det\SQ_{\{x\}\times
T,i})^{d_i(x)}\rbrace\otimes \bigotimes^2_{j=1}
(det\Cal F_{\{y_j\}\times T})^{\ell_j}$$
where $\pi_T$ is the projection $X_0\times T\to T$, and
$det\,R\pi_T\Cal F$
is the determinant bundle defined as
$$\{det\,R\pi_T\Cal
F\}_t:=\{det\,H^0(X,\Cal F_t)\}^{-1}\otimes\{
det\,H^1(X,\Cal F_t)\}.$$

\proclaim{Theorem 1.2} There is an unique ample line bundle
$\Theta_{\Cal U_{X_0}}=\Theta(k,\ell_1,\ell_2,\vec a,\vec n,\alpha,I)$ 
on $\Cal U_{X_0}$ such that for any given family of semistable parabolic 
sheaf $\Cal F$ parametrised by $T$, we have
$\phi_T^*\Theta_{\Cal U_{X_0}}=\Theta_{\SF}$, where $\phi_T$ is the 
induced map $T\to\Cal U_{X_0}$\endproclaim

\demo{Proof} By using the descendant lemma (see the next Lemma 1.2), we will
show that the line bundle $\Theta_{\SR^{ss}}:=\Theta_{\SE}$ on $\SR^{ss}$ 
descends to the required ample line $\Theta_{\Cal U_{X_0}}$, where $\SE$ is
a universal quotient over $X_0\times\SR^{ss}$.

We known that the stablizer $stab(q)=\lambda\cdot id$ for $q\in\SR^{s}$, which acts on
$\Theta_{\SR^{ss}}$ via
 $$\lambda^{-k\chi+ \sum_{x\in
I}\sum^{l_x}_{i=1}d_i(x)r_i(x)+r\sum_{x\in
I}
\alpha_x+r(\ell_1+\ell_2)}=\lambda^0=1.$$
If $q\in\SR^{ss}\setminus\SR^s$ with closed orbit, we known that
$$\SE_q=m_1E_1\oplus m_2E_2\oplus\cdots\oplus m_tE_t$$
with $par\mu_m(E_j)=par\mu_m(\SE_q)$, which means that (if assuming $E_j$
of rank $(r_1,r_2)$) 
$$-k\chi(E_j)+r_1\sum_{x\in I_1}\alpha_x+r_2\sum_{x\in I_2}\alpha_x
+\sum_{x\in I}\sum^{l_x}_{i=1}d_i(x)dim\frac{E_{j\,x}}{E_{j\,x}\cap F_i(E)_x}
+r_1\ell_1+r_2\ell_2=0.$$
Thus $(\lambda_1\cdot id_{m_1},\cdots,\lambda_t\cdot id_{m_t})\in 
stab(q)=GL(m_1)\times\cdots\times GL(m_t)$ acts trivially on $\Theta_{\SR^{ss}}$, which implies that the $stab(q)$ acts trivially on
$\Theta_{\SR^{ss}}$ and thus descends to a line bundle $\Theta_{\SU_{X_0}}$
having the required universal property.

To show the ampleness of $\Theta_{\SU_{X_0}}$, noting that 
$detR\pi_{\SR^{ss}}\Cal E(N)$ is trivial and
$$detR\pi_{\SR^{ss}}\Cal E=(det\SE_{y_1})^{c_1N}\otimes(det\SE_{y_2})^{c_2N}
\otimes detR\pi_{\SR^{ss}}\Cal E(N),$$
we see that the restriction of the polarization to $\SR^{ss}$ is
$$(detR\pi_{\SR^{ss}}\Cal E(m))^{\frac{\ell}{m-N}}\otimes\bigotimes_{x\in I}
\left\{(det\Cal E_x)^{\alpha_x}\otimes\bigotimes^{l_x}_{i=1}
(det\SQ_x)^{d_i(x)}\right\}=\Theta_{\SR^{ss}}.$$
Thus, by general theorems of GIT, some power of $\Theta_{\SR^{ss}}$ descends
to an ample line bundle, which implies that some power of $\Theta_{\SU_{X_0}}$
is ample. 
\enddemo

\proclaim{Lemma 1.2} Let $G$ be a reductive algebraic group and $V$ a scheme
with $G$-action. Suppose that there exists the good quotient
$\pi:V\to V//G$. Then a vector bundle $E$ with $G$-action over $V$ descends
to $V//G$ iff stablizer $stab(y)$ of $y$ acts on $E_y$ trivially for any $y\in V$ with closed
orbit.\endproclaim

It is known that for any torsion free sheaf $F$  of rank $(r_1,r_2)$ on
$X_0$ there are integers $a$, $b$, $c$ such that 
$$F_{x_0}\cong \SO^a_{X_0,x_0}\oplus\SO^b_{X_1,x_0}\oplus\SO^c_{X_2,x_0},$$
where $a$, $b$, $c$ are determined uniquely and satisfying
$$r_1=a+b,\quad r_2=a+c,\quad dim(F_{x_0}\otimes k(x_0))=a+b+c.$$
Thus we can define $\bold a(F):=a$ for any torsion free sheaf $F$ on $X_0$, and
we have

\proclaim{Lemma 1.3} Let $0\to G\to
F\to E\to 0$ be an exact
sequence of torsion free sheaves on $X_0$. Then
$$\bold
a(F)\ge\bold a(G)+\bold a(E).$$
\endproclaim

\demo{Proof} It is clear by counting the dimension of their fibres at $x_0$.
\enddemo 
 
Let
$\SR_a=\{F\in\SR|\, F\otimes\hat\SO_{x_0}=
\hat\SO_{x_0}^{\oplus a}\oplus
m_{x_0}^{\oplus(r-a)}\}$,
and
$\hat\SW_i=\SR_0\cup\SR_1\cup\cdots\cup\SR_i$
(which are closed in
$\SR$) endowed with their reduced scheme structures.
The subschemes
$\hat\SW_i$ are $SL(n)$-invariant, and yield closed reduced
subschemes of
$\SU_X$. It is clear
that
$$\SR\supset\hat\SW_{r-1}\supset\hat\SW_{r-2}\supset\cdots\supset\hat\SW_1
\supset\hat\SW_0=\SR_0$$
$$\SU_X\supset\SW_{r-1}\supset\SW_{r-2}\supset\cdots\supset\SW_1
\supset\SW_0.$$

Let $q_0\in\SR$ be a point corresponding to a torsion free sheaf $\SF_0$ such
that 
$$\SF_0\otimes\SO_{X_0,x_0}\cong m_{x_0}^{r-a_0}\oplus\SO_{X_0,x_0}^{a_0}.$$
We consider the variety
$$Z=\{(X,Y)\in M(r-a_0)\times M(r-a_0)\,|\,X\cdot Y=Y\cdot X=0\},$$
and its subvarieties $Z'=\{(X,Y)\in Z\,|\,rk(X)+rk(Y)\le a\}.$ 
Then the reduced coordinate ring of $Z$ is 
$$\Bbb C[Z]:=\frac{\Bbb C[X, Y]}{(XY, YX)},$$
where $X:=(x_{ij})_{r-a_0\times r-a_0}$ and $Y:=(y_{ij})_{r-a_0\times r-a_0}$
(see Lemma 4.8 of [Su]), and $Z'$ is a union of reduced subvarieties 
of $Z$ (see the proof of Theorem 4.2 in [Su]). Thus we can sum up the
the arguments of [NS] and [Su] (see also [Fa]) into a lemma

\proclaim{Lemma 1.4} The variety $Z$, $Z'$ is respectively the local model
of $R,$ $\hat\SW_a$ at the point $q_0$. More precisely, there are some
integers $s$ and $t$ such that
$$\hat\SO_{\SR,q_0}[[u_1,\cdots,u_s]]\cong \hat\SO_{Z,(0,0)}
[[v_1,\cdots,v_t]],$$
$$\hat\SO_{\hat\SW_a,q_0}[[u_1,\cdots,u_s]]\cong \hat\SO_{Z',(0,0)}
[[v_1,\cdots,v_t]].$$
In particular, $\SW_a$ ($0\le a\le r$) are reduced and seminormal.\endproclaim

\heading \S 2 Moduli space of generalized parabolic sheaves \endheading

Let $\pi:\widetilde X_0\to X_0$ be the normalisation of $X_0$ and
$\pi^{-1}
(x_0)=\{x_1,x_2\},$ then $\widetilde X_0$ is a disjoint union of
$X_1$ and $X_2$, any coherent sheaf $E$ on $\widetilde X_0$ is determined
by a pair $(E_1,E_2)$ of coherent sheaves on $X_1$ and $X_2$. We call as
before that $E$ is of rank $(r_1,r_2)$ if $E_i$ has rank $r_i$ on $X_i$ 
($i=1,2$) and define the rank of $E$ to be 
$$r(E):=\frac{c_1r_1+c_2r_2}{c_1+c_2}.$$
We can also define similarly the modified parabolic Euler characteristic
$par\chi_m(E)$ if $E$ has parabolic structures at points $x\in\pi^{-1}(I)$
(we will identify $I$ with $\pi^{-1}(I)$, and note that $m(E)$ defined
in Definition 1.3 is only depend on $r_1$ and $r_2$ since $\SO(1)$, $\alpha$
and $\vec a(x)$ are fixed). 
 
\proclaim{Definition 2.1} A generalized parabolic sheaf of rank $(r_1,r_2)$
(abbreviated to GPS) 
$$\bold E:=(E, E_{x_1}\oplus E_{x_2}@>q>>Q)$$
on $\widetilde X_0$ is a coherent sheaf $E$ on $\widetilde X_0$, torsion free
of rank $(r_1,r_2)$ outside $\{x_1,x_2\}$ with parabolic structures at points
$\{x\}_{x\in I}$, together with a quotient $E_{x_1}\oplus E_{x_2}@>q>>Q$.
A morphism $f:(E, E_{x_1}\oplus E_{x_2}@>q>>Q)\to
(E', E'_{x_1}\oplus E'_{x_2}@>q'>>Q')$ of GPS
is a morphism $f: E\to E'$ of parabolic sheaves, which maps $ker(q)$ into
$ker(q')$.\endproclaim

We will consider the generalized parabolic sheaves $(E, Q)$
of rank $r_1=r_2=r$ and $dim(Q)=r$ with parabolic structures 
of type $\{\vec n(x)\}_{x\in I}$ and weights
$\{\vec a(x)\}_{x\in I}$ at the points of $\pi^{-1}(I)$, and we will call them 
the GPS of rank $r$. Furthermore, by a family of GPS of rank $r$ over
$T$, we
mean the following
\roster
\item a rank $r$ sheaf $\Cal E$ on
$\widetilde X_0\times T$ flat over
$T$ and locally free outside
$\{x_1,x_2\}\times T$.
\item a locally free rank $r$ quotient $\Cal Q$ of
$\Cal E_{x_1}\oplus
\Cal E_{x_2}$ on $T$.
\item a flag bundle $Flag_{\vec n(x)}(\Cal
E_x)$ on $T$ with given weights for each
$x\in
I$.\endroster

\proclaim{Definition 2.2} A GPS $(E,Q)$ is called
semistable
(resp., stable), if for every nontrivial subsheaf $E'\subset E$
such that $E/E'$ is torsion
free outside $\{x_1,x_2\},$ we
have, with the induced parabolic structures at points $\{x\}_{x\in I}$,
$$par\chi_m(E')-dim(Q^{E'})\leq
rk(E')\cdot\frac{par\chi_m(E)-dim(Q)}{rk(E)}
\,\quad (\text{resp.,
$<$}),$$
where $Q^{E'}=q(E'_{x_1}\oplus E'_{x_2})\subset
Q.$ \endproclaim

Let $\chi_1$ and $\chi_2$ be integers such that $\chi_1+\chi_2-r=\chi$, and fix,
for $i=1,2$, the polynomials $P_i(m)=c_irm+\chi_i$ and $\SW_i=\SO_{X_i}(-N)$
where $\SO_{X_i}(1)=\SO (1)|_{X_i}=\SO_{X_i}(c_iy_i)$. 
Write $V_i=\Bbb C^{P_i(N)}$
and consider the Quot schemes $Quot(V_i\otimes\SW_i, P_i)$, 
let $\bwq_i$ be the
closure of the open set 
$$\bold Q_i=\left\{\aligned&\text{$V_i\otimes\SW_i\to E_i\to 0$, with 
$H^1(E_i(N))=0$ and}\\&\text{$V\to H^0(E_i(N))$ 
induces an isomorphism}\endaligned\right\},$$
we have the universal quotient $V_i\otimes\SW_i\to \SF^i\to 0$ on 
$X_i\times\bwq_i$ and the relative flag scheme
$$\SR_i=\underset x\in I_i \to{\times_{\bwq_i}}
Flag_{\vec n(x)}(\SF^i_x)\to\bwq_i.$$
Let $\SE^i$ be the pullback of $\SF^i$ to $X_i\times\SR_i$ and 
$$\rho:\widetilde\SR=Grass_r(\SE^1_{x_1}\oplus\SE^2_{x_2})\to\SR_1\times\SR_2.$$
Then we see that, for $N$ large enough, every semistable GPS appears as a point
of $\widetilde\SR$. To rewrite $\SR_1\times\SR_2$ so that it unified
the $R$ in last section, let $V=V_1\oplus V_2$, $\SF=\SF^1\oplus\SF^2$ and
$\SE=\SE^1\oplus\SE^2$, we have
$$\SR_1\times\SR_2= \underset x\in I \to{\times_{\bwq_1\times\bwq_2}}
Flag_{\vec n(x)}(\SF_x)\to\bwq_1\times\bwq_2.\tag2.1$$
Note that $V_1\otimes\SW_1\oplus V_2\otimes\SW_2\to\SF\to 0$ is a 
$\bwq_1\times\bwq_2$-flat quotient with Hilbert polynomial $P(m)=P_1(m)+P_2(m)$
on $\widetilde X_0\times (\bwq_1\times\bwq_2)$, we have for $m$ large enough
a $G$-equivariant embedding
$$\bwq_1\times\bwq_2\hookrightarrow Grass_{P(m)}(V_1\otimes W_1^m
\oplus V_2\otimes W_2^m),$$
where $W_i^m=H^0(\SW_i(m))$ and $G=(GL(V_1)\times GL(V_2))\cap SL(V).$ 

A (closed) point $(p=p_1\oplus p_2,\{p_{r(x)},p_{r_1(x)},...,p_{r_{l_x}(x)}\}_{x\in I})$ of $\SR_1\times\SR_2$ by the expression of (2.1)
is given by points $V_i\otimes\SW_i@>p_i>>E^i\to 0$ of the Quot schemes ($i=1,2$), together with
quotients (if we write $\Cal V_{\widetilde X_0}=V_1\otimes\SW_1\oplus V_2\otimes\SW_2$ and $E=E^1\oplus E^2$) 
$$\{\Cal V_{\widetilde X_0}@>p_{r(x)}>>Q_{r(x)},
\Cal V_{\widetilde X_0}@>p_{r_1(x)}>>Q_{r_1(x)},...,
\Cal V_{\widetilde X_0}@>p_{r_{l_x(x)}}>>Q_{r_{l_x}(x)}\}_{x\in I},$$
where $r_i(x)=dim(E_x/F_i(E)_x)=n_1(x)+\cdots+n_i(x)$, and $Q_{r(x)}:=E_x,$
$Q_{r_1(x)}:=E_x/F_1(E)_x, ..., Q_{r_{l_x}(x)}:=E_x/F_{l_x}(E)_x$, the morphisms
$p_{r(x)}$ and $p_{r_j(x)}$ ($j=1,...,l_x$) are defined to be
$$p_{r(x)}:\Cal V_{\widetilde X_0}@>p>>E\to E_x,
\quad p_{r_j(x)}:\Cal V_{\widetilde X_0}@>p_{r(x)}>>Q_{r(x)}=E_x\to E_x/F_j(E)_x.$$
Thus
we have a $G$-equivariant embedding 
$$\SR_1\times\SR_2\hookrightarrow Grass_{P(m)}(V_1\otimes W_1^m
\oplus V_2\otimes W_2^m)\times\bold {Flag},$$
where $\bold{Flag}$ is defined to be
$$\bold{Flag}=\prod_{x\in I}\{Grass_{r(x)}(V)\times Grass_{r_1(x)}(V)
\times\cdots\times
Grass_{r_{l_x}(x)}(V)\},$$
which maps a point $(p=p_1\oplus p_2,\{p_{r(x)},p_{r_1(x)},...,p_{r_{l_x}(x)}\}_{x\in I})$
of $\SR_1\times\SR_2$ to the point
$$(H^0(\Cal V_{\widetilde X_0}(m))@>g>>U,\{V@>g_{r(x)}>>U_{r(x)},
V@>g_{r_1(x)}>>U_{r_1(x)},...,
V@>g_{r_{l_x}(x)}>>U_{r_{l_x}(x)}\}_{x\in I})$$
of $Grass_{P(m)}(V_1\otimes W_1^m\oplus
V_2\otimes W_2^m)\times\bold {Flag}$, 
where $g:=H^0(p(m))$, $U:=H^0(E(m))$, $g_{r(x)}:=H^0(p_{r(x)}(N)),$ $U_{r(x)}:=H^0(Q_{r(x)}),$ and $g_{r_j(x)}:=H^0(p_{r_j(x)}(N))$, $U_{r_j(x)}:=H^0(Q_{r_j(x)})$ ($j=1,...,l_x$). Finally, we get a 
$G$-equivariant embedding 
$$\widetilde\SR\hookrightarrow\bold G'= Grass_{P(m)}(V_1\otimes W_1^m
\oplus V_2\otimes W_2^m)\times\bold {Flag}\times Grass_r(V_1\oplus V_2)$$
as follows: a point of $\widetilde\SR$ is given by a point of $\SR_1\times\SR_2$
together with a quotient $E_{x_1}\oplus E_{x_2}@>q>>Q$, then above embedding
maps $E_{x_1}\oplus E_{x_2}@>q>>Q$ to
$$g_G:=H^0(q(N)): V_1\oplus V_2=H^0(\Cal V_{\widetilde X_0}(N))\to H^0(E(N))
\to E_{x_1}\oplus E_{x_2}@>q>> Q.$$
Given $\bold G'$ the polarisation (using
the obvious notation):
$$\left\{\frac{\ell}{m-N}\times\prod_{x\in
I}\{\alpha_x, d_1(x),\cdots,
d_{l_x}(x)\right\}\times k,$$
we have the analogue of Proposition 1.1, whose proof (we refer to Proposition 1.14 of [B3], or Lemma 5.4 of [NS]) is a modification of Theorem 4.17 in [Ne]
since our group $G$ here is different from that of [Ne].
 
\proclaim{Proposition 2.1} A point 
$(g,\{g_{r(x)},g_{r_1(x)},...,g_{r_{l_x}(x)}\}_{x\in I},g_G)\in \bold G'$ is stable
(respectively, semistable) for the action of $G$, with respect to the 
above polarisation (we refer to this from now on as GIT-stability), iff for all
nontrivial subspaces $H\subset V$, where $H=H_1\oplus H_2$ and $H_i\subset V_i$
($i=1,2$), we have (with $h=dim H$ and $\bar H:=H_1\otimes W^m_1\oplus H_2\otimes W_2^m$)
$$\aligned &\frac{\ell}{m-N}\left(hP(m)-P(N)dim\,g(\bar H)\right)+
\sum_{x\in I}\alpha_x(rh-P(N)dim\,g_{r(x)}(H))+\\
&\sum_{x\in I}\sum^{l_x}_{i=1}d_i(x)\left(r_i(x)h-P(N)dim\,g_{r_i(x)}(H)\right)
+k\left(rh-P(N)dim\,g_G(H)\right)<(\le)\,0.\endaligned$$\endproclaim

\proclaim{Proposition 2.2} Suppose 
$(p,\{p_{r(x)},p_{r_1(x)},...,p_{r_{l_x}(x)}\}_{x\in I},q)\in\widetilde\SR$ 
is a point such that $E$ is torsion free outside $\{x_1,x_2\}$. Then 
$\bold E=(E,E_{x_1}\oplus E_{x_2}@>q>>Q)$ is stable (respectively, semistable)
iff for every subsheaf $0\neq F\neq E$ we have (using the notation 1.2)
$$\aligned \frac{\ell}{m-N}(&\chi(F(N))P(m)-P(N)\chi(F(m)))+
\sum_{x\in I}\alpha_x(r\chi(F(N))-P(N)h^0(Q_{r(x)}^F))\\
&+\sum_{x\in I}\sum^{l_x}_{i=1}d_i(x)(r_i(x)\chi(F(N))-P(N)h^0(Q_{r_i(x)}^F))
\\&+k(r\chi(F(N))-P(N)dim(Q^F))<(\le)\,0.\endaligned$$\endproclaim

\demo{Proof} For subsheaf $F\subset E$ such that $E/F$ is torsion free outside
$\{x_1,x_2\}$, by the same computation in Proposition 1.2, we have
$$LHS(F)=kP(N)\left(par\chi_m(F)-dim(Q^F)-r(F)\frac{par\chi_m(E)-r}{r}\right).$$
Thus $\bold E$ is stable (semistable) iff $LHS(F)<(\le)0$ for the required $F.$
If $E/F$ has torsion outside $\{x_1,x_2\}$, then $LHS(F)<0.$\enddemo

\proclaim{Lemma 2.1} There exist $N$ and $M_1(N)$ such that for $m\ge M_1(N)$ the
following holds. Suppose $(p,\{p_{r(x)},p_{r_1(x)},...,p_{r_{l_x}(x)}\}_{x\in I}, q)\in\widetilde\SR$ is a point which is GIT-semistable then for all quotients 
$E@>T>>\Cal G\to 0$ we have (with $Q^{\Cal G}:=Q/q(ker(T))$)
$$\aligned h^0(\Cal G(N))\ge&\frac{1}{k}\left((c_1+c_2)r(\Cal G)\ell+
\sum_{x\in I}\alpha_xh^0(Q^{\Cal G}_{r(x)})+\sum_{x\in I}\sum^{l_x}_{i=1}
d_i(x)h^0(Q^{\Cal G}_{r_i(x)})\right)\\&+h^0(Q^{\Cal G}).\endaligned$$ 
In particular, $E$ is torsion free outside $\{x_1,x_2\}$, $q$ maps the torsion
on $\{x_1,x_2\}$ to $Q$ injectively and
$V\to H^0(E(N))$ is an isomorphism. \endproclaim

\demo{Proof} The proof of Lemma 1.1 goes through with obvious modifications
except that we can not assume that the sheaves $E$ are torsion free at
$x_1$ and $x_2.$ To see it clearly, we write out the proof of $E$ being 
torsion free outside $\{x_1,x_2\}.$

Let $\tau=Tor(E)$, $\Cal G=E/\tau$ and applying the above inequality, noting that
$h^0(\Cal G(N))=P(N)-h^0(\tau)$, $h^0(Q^{\Cal G}_{r(x)})=r-h^0(Q^{\tau}_{r(x)})$ and
$h^0(Q^{\Cal G}_{r_i(x)})=r_i(x)-h^0(Q^{\tau}_{r_i(x)}),$ we have
$$kh^0(\tau)\le k\cdot dim(Q^{\tau})+\sum_{x\in I}(\alpha_x+a_{l_x+1}(x)-a_1(x))h^0(\tau_x),$$
by which one can conclude that $\tau=0$ outside $\{x_1,x_2\}$ and $h^0(\tau_{x_1}\oplus\tau_{x_2})-dim(Q^{\tau})=0$ since $\alpha_x<k-a_{l_x+1}(x)+a_1(x).$ In particular, $q$ maps the torsion on 
$\{x_1,x_2\}$ to $Q$ injectively.\enddemo

\remark{Remark 2.1} The proof of Lemma 1.1 and Lemma 2.1 actually implies
that one can take $N$ big enough such that for a GIT-semistable point the
sheaf $E$ involved satisfies the condition $H^1(E(N)(-x-x_1-x_2))=0$ for
any $x\in X_0$, which implies that $E(N)$ and $E(N)(-x_1-x_2)$ are generated 
by global sections and $H^0(E(N))\to E(N)_{x_1}\oplus E(N)_{x_2}$ 
is surjective. Conversely, it is easy to prove that every semistable GPS will 
satisfy above conditions if $N$ big enough.\endremark

\proclaim{Proposition 2.3} There exist integers $N>0$ and $M(N)>0$ such that
for $m\ge M(N)$ the following is true. A point
$(p,\{p_{r(x)},p_{r_1(x)},...,p_{r_{l_x}(x)}\}_{x\in I},q)\in\widetilde\SR$ 
is GIT-stable (respectively, GIT-semistable) iff 
the quotient $E$ is torsion free outside $\{x_1,x_2\}$ and $\bold E=(E,q)$
is stable (respectively, semistable) GPS, the map $V\to H^0(E(N))$ is an isomorphism.\endproclaim

\demo{Proof} The proof is the same with that of Proposition 1.3 by some
obvious notation modifications.\enddemo

\proclaim{Notation 2.1} Define $\SH$ to be the subscheme of
$\wt\SR$
parametrising the generalised parabolic sheaves $\bold E=(E,E_{x_1}
\oplus E_{x_2}@>q>>Q)$ satisfying
\roster\item $\Bbb C^{P(N)}\cong H^0(E(N)),$
and $H^1(E(N)(-x_1-x_2-x))=0$
for any $x\in\wt X_0$
\item $\text{Tor}E$ is
supported on $\{x_1,x_2\}$ and
$(\text{Tor}E)_{x_1}\oplus(\text{Tor}E)_{x_2}\hookrightarrow
Q.$
\endroster
Let $\wt\SR^{ss}$ ($\wt\SR^s$) be the open set of $\wt\SR$
consists the semistable (stable) GPS, then it is clear that
$$\wt\SR^{ss}\overset\text{open}\to\hookrightarrow
\SH\overset\text{open}\to\hookrightarrow
\wt\SR.$$
\endproclaim

We will introduce the so called $s$-equivalence of GPS later in Definition 2.6.
It is also known that $\SH$ is reduced, normal
and Gorenstein with only rational singularities (see Proposition 3.2 and
Remark 3.1 in [Su]).

\proclaim{Theorem 2.1} For given datas in Notation 1.1 satisfying $(*)$ and
$\chi_1$, $\chi_2$ with $\chi_1+\chi_2-r=\chi$, there exists an irreducible,
Gorenstein, normal projective variety $\SP_{\chi_1,\chi_2}$ with only rational
singularities, which is the coarse moduli space of $s$-equivalence classes of
semistable GPS $(E,Q)$ on $\wt X_0$ of rank $r$ and $\chi(E_j)=\chi_j$ 
($j=1,2$) with parabolic structures of type $\{\vec n(x)\}_{x\in I}$ and
weights $\{\vec a(x)\}_{x\in I}$ at points $\{x\}_{x\in I}$.\endproclaim

\demo{Proof}  The existence of the moduli space
and its projectivity follows above Proposition 2.3 and G.I.T., the
other properties follow the corresponding properties of $\SH$ and the fact
that $\wt\SR^{ss}\subset\SH$ for if $N$ big enough.\enddemo

Recall that we have the universal quotient $\SE^1$ on $X_1\times\SR_1$, flat
over $\SR_1$, and torsion free of rank $r$ outside $\{x_1\}$ with Euler
characteristic $\chi_1$, together with, for each $x\in I_1$, a flag
$$\aligned\SE^1_{\{x\}\times\SR_1}=F_0(\SE^1_{\{x\}\times\SR_1})\supset
F_1(\SE^1_{\{x\}\times\SR_1})\supset\cdots\supset 
F_{l_x}(\SE^1_{\{x\}\times\SR_1})\supset &F_{l_x+1}(\SE^1_{\{x\}\times\SR_1})\\&=0\endaligned$$
of subbundles of type $\vec n(x)$ and weights $\vec a(x)$. Let
$\SQ_{x,i}=\SE^1_{\{x\}\times\SR_1}/F_i(\SE^1_{\{x\}\times\SR_1})$, we can 
define a line bundle $\Theta_{\SR_1}$ on $\SR_1$ as
$$(detR\pi_{\SR_1}\SE^1)^k\otimes\bigotimes_{x\in I_1}\left\{(det
\SE^1_{\{x\}\times\SR_1})^{\alpha_x}\otimes\bigotimes^{l_x}_{i=1}
(det\SQ_{x,i})^{d_i(x)}\right\}\otimes(det\SE^1_{\{y_1\}
\times\SR_1})^{\ell_1}.$$
Similarly, we can define the line bundle $\Theta_{\SR_2}$ on $\SR_2$ and the
$G$-line bundle
$$\Theta_{\wt\SR}:=\rho^*(\Theta_{\SR_1}\otimes\Theta_{\SR_2})
\otimes(det\SQ)^k$$
on $\wt\SR$, where $\rho^*(\SE^1_{x_1}\oplus\SE^2_{x_2})\to\SQ\to 0$ is the
universal quotient on $\wt\SR$. One can check that $\Theta_{\wt\SR}$ is the
restriction of ample polarisation used to linearize the action of $G$, thus
some power of $\Theta_{\wt\SR}$ descends to an ample line bundle on $\SP_{\chi_1,\chi_2}$. In fact, we have 

\proclaim{Lemma 2.2} The $\Theta_{\wt\SR^{ss}}$ descends to an ample line
bundle $\Theta_{\SP_{\chi_1,\chi_2}}$ on $\SP_{\chi_1,\chi_2}$.\endproclaim 

\demo{Proof} The proof is similar with Theorem 1.2, we only make a remark
here. If $(E,Q)$ is a semistable GPS of rank $r$ and $(E',Q')$ a sub-GPS
of $(E,Q)$ with 
$$par\chi_m(E')-dim(Q')=r(E')\cdot\frac{par\chi_m(E)-dim(Q)}{r},$$
we have (assuming that $E'$ is of rank $(r_1,r_2)$)
$$\aligned&
-k\chi(E')+r_1\sum_{x\in I_1}\alpha_x+r_2\sum_{x\in I_2}\alpha_x
+\sum_{x\in I}\sum^{l_x}_{i=1}d_i(x)dim\frac{E'_x}{E'_x\cap F_i(E_x)}
+r_1\ell_1+r_2\ell_2\\
&+k\cdot dim(Q') 
\\&=\frac{-k\chi+r\sum_{x\in I}\alpha_x+\sum_{x\in I}\sum^{l_x}_{i=1}d_i(x)r_i(x)+r(\ell_1+\ell_2)}
{r}\cdot r(E')
=0.\endaligned$$
\enddemo

\proclaim{Notation 2.2} Let $\SR_{1\,F}\subset\SR_1$ ($\SR_{2\,F}\subset\SR_2$)
be the open set of points corresponding the vector bundles on $X_1$ ($X_2$),
and $\wt\SR_F={\rho}^{-1}(\SR_{1\,F}\times\SR_{2\,F}),$ then
$$\rho:\wt\SR_F\to\SR_{1\,F}\times\SR_{2\,F}$$
is a grassmannian bundle over $\SR_{1\,F}\times\SR_{2\,F}$,
and $\wt\SR_F\subset\SH.$ We
define
$$R^1_{F,a}:=\{(E,Q)\in\wt\SR_F|\text{$E_{x_1}\to Q$ has rank
$a$}\},$$
and $\hat\SD_{F,1}(i):=R^1_{F,0}\cup\cdots\cup R^1_{F,i},$ which
have the
natural scheme structures. The subschemes $R^2_{F,a}$
and
$\hat\SD_{F,2}(i)$ are defined similarly. Let $\hat\SD_1(i)$ and
$\hat\SD_2(i)$ be the zariski closure of $\hat\SD_{F,1}(i)$
and
$\hat\SD_{F,2}(i)$ in $\wt\SR$. Then they are reduced, irreducible
and
$G$-invariant closed subschemes of $\wt\SR$, thus
inducing
closed subschemes $\SD_1(i)_{\chi_1,\chi_2}$, $\SD_2(i)_{\chi_1,\chi_2}$ of $\SP_{\chi_1,\chi_2}$. Clearly, we
have 
(for $j=1,2$)
that
$$\wt\SR\supset\hat\SD_j(r-1)\supset\hat\SD_j(r-2)\cdots\supset
\SD_j(1)\supset\hat\SD_j(0)$$
$$\SP_{\chi_1,\chi_2}\supset\SD_j(r-1)_{\chi_1,\chi_2}\supset\SD_j(r-2)_{\chi_1,\chi_2}\supset\cdots\SD_j(1)_{\chi_1,\chi_2}
\supset\SD_j(0)_{\chi_1,\chi_2}.$$
\endproclaim

\proclaim{Lemma 2.3} $\SH$, $\hat\SD_j(a)$ and $\hat\SD_1(a)\cap\hat\SD_2(b)$ are
reduced, normal with rational singularities. In particular,
$\SP_{\chi_1,\chi_2}$, $\SD_j(a)_{\chi_1,\chi_2}$ and $\SD_1(a)_{\chi_1,\chi_2}\cap\SD_2(b)_{\chi_1,\chi_2}$ are reduced, normal 
with rational singularities.\endproclaim

\demo{Proof} This is the copy of Proposition 3.2 in [Su] and the proof there
goes through.\enddemo

Let $(E,Q)$ be a semistable GPS of rank $r$ with $E=(E_1, E_2)$ and
$\chi_j=\chi(E_j)$ ($j=1,2$). Then, by the definition of semistability,
we have (for $j=1,2$) that
$$par\chi_m(E_j)-dim(Q^{E_j})\le \frac{c_j}{c_1+c_2}(par\chi_m(E)-r).$$
Recall that $\chi_1+\chi_2-r=\chi$ and $n_j$ ($j=1,2$) denotes
$$\frac{1}{k}\left(\sum_{x\in I_j}\sum^{l_x}_{i=1}
d_i(x)r_i(x)+r\sum_{x\in I_j}\alpha_x+r\ell_j\right),$$
we can rewrite the above inequality into
$$\aligned&n_1+r-dim(Q^{E_2})\le\chi(E_1)\le n_1+dim(Q^{E_1})\\
&n_2+r-dim(Q^{E_1})\le\chi(E_2)\le n_2+dim(Q^{E_2}).\endaligned\tag2.1$$
Thus, for fixed $\chi$, the moduli space of $s$-equivalence classes of 
semistable GPS $(E,Q)$ on $\wt X_0$ of rank $r$ and $\chi(E)=\chi+r$ with
parabolic structures of type $\{\vec n(x)\}_{x\in I}$ and weights
$\{\vec a(x)\}_{x\in I}$ at points $\{x\}_{x\in I}$ is the disjoint union
$$\SP:= \coprod_{\chi_1+\chi_2=\chi+r}\SP_{\chi_1,\chi_2},$$
where $\chi_1$, $\chi_2$ satisfy the inequalities
$$n_1\le\chi(E_1)\le n_1+r,\quad n_2\le\chi(E_2)\le n_2+r.\tag2.2$$

\proclaim{Notation 2.3} The ample line bundles $\{\Theta_{\SP_{\chi_1,\chi_2}}\}$ determine an ample line bundle $\Theta_{\SP}$
on $\SP$, and for any $0\le a\le r$, we define the subschemes 
$$\SD_1(a):= \coprod_{\chi_1+\chi_2=\chi+r}\SD_1(a)_{\chi_1,\chi_2},
\quad \SD_2(a):= \coprod_{\chi_1+\chi_2=\chi+r}\SD_2(a)_{\chi_1,\chi_2}.$$
We will simply write $\SD_1:=\SD_1(r-1)$ and $\SD_2:=\SD_2(r-1)$.\endproclaim

In order to introduce a sheaf theoretic description of the so called $s$-equivalence of GPS, we enlarge the category by considering all of the
GPS including the case $r(E)=0$, and also assume that $|I|=0$ for simplicity.

\proclaim{Definition 2.3} A GPS $(E,Q)$ is called
semistable
(resp., stable), if \roster
\item when $rank(E)>0$, then for every nontrivial subsheaf $E'\subset E$
such that $E/E'$ is torsion
free outside $\{x_1,x_2\},$ we
have, with the induced parabolic structures at points $\{x\}_{x\in I}$,
$$par\chi_m(E')-dim(Q^{E'})\leq
rk(E')\cdot\frac{par\chi_m(E)-dim(Q)}{rk(E)}
\,\quad (\text{resp.,
$<$}),$$
where $Q^{E'}=q(E'_{x_1}\oplus E'_{x_2})\subset
Q.$
\item when $rank(E)=0$, then $E_{x_1}\oplus E_{x_2}=Q$ (resp. $E_{x_1}\oplus E_{x_2}=Q$ and $dim(Q)=1$)\endroster 
\endproclaim

\proclaim{Definition 2.4} If $(E,Q)$ is
a GPS and $rank(E)>0$, we
set
$$\mu_G[(E,Q)]=\frac{deg(E)-dim(Q)}{rank(E)}.$$\endproclaim

It is
useful to think of an $m$-GPS as a sheaf $E$ on $\wt X_0$ together
with a map
$\pi_*E\to\,_{x_0}Q\to 0$ and $h^0(_{x_0}Q)=m$. Let $K_E$
denote the kernel
of $\pi_*E\to Q$.

\proclaim{Definition 2.5}
Given an exact sequence
$$0\to E'\to E\to E''\to 0$$
of sheaves on $\wt X$,
and $\pi_*E\to Q\to 0$ a generalised parabolic
structure on $E$, we define
the generalised parabolic structures on $E'$
and $E''$ via the
diagram
$$\CD
0@>>> \pi_*E'@>>>\pi_*E@>>> \pi_*E''@>>>  0  \\
@.    @VVV
@VVV        @VVV        @. \\
0@>>> Q'    @>>>   Q   @>>>    Q''  @>>>
0
\endCD$$    
The first horizontal sequence is exact because $\pi$ is
finite, $Q'$
is defined as the image in $Q$ of $\pi_*E'$ so that the first
vertical
arrow is onto, $Q''$ is defined by demanding that the second
horizontal
sequence is exact, and finally the third vertical arrow is onto
by the
snake lemma. We will write
$$0\to (E',Q')\to (E,Q)\to (E'',Q'')\to
0$$
whose meaning is clear.
\endproclaim

\proclaim{Proposition 2.4} Fix a
rational number $\mu$. Then the category
$\Cal C_{\mu}$ of semistable GPS
$(E,Q)$ such that $rank(E)=0$ or,
$rank(E)>0$ with $\mu_G[(E,Q)]=\mu$, is
an abelian, artinian, noetherian
category whose simple objects are the
stable GPS in the
category.\endproclaim

One can conclude, as usual, that
given a semistable GPS $(E,Q)$ it has a
Jordan-Holder filtration, and the
associated graded GPS $gr(E,Q)$ is
uniquely determined by $(E,Q)$. Thus we
have

\proclaim{Definition 2.6} Two semistable GPS $(E_1,Q_1)$ and
$(E_2,Q_2)$
are said to be $s$-equivalent if they have the same associated
graded GPS,
namely,
$$(E_1,Q_1)\sim  (E_2,Q_2)\quad\Longleftrightarrow\quad
gr(E_1,Q_1)\cong
gr(E_2,Q_2).$$\endproclaim

\remark{Remark 2.2} Any stable
GPS $(E,Q)$ with $rank(E)>0$ must
be locally free (i.e., E is locally free), and two stable GPS are s-equivalent
iff they are isomorphic.\endremark

\proclaim{Proposition 2.5} Every semistable $(E',Q')$ with $rank(E')>0$
is $s$-equivalent to a semistable $(E,Q)$ with $E$ locally free. Moreover,
\roster\item if $E'$ has torsion of dimension $t$ at $x_2$, then $(E',Q')$
is $s$-equivalent to a semistable $(E,Q)$ with $E$ locally free and
$$rank(E_{x_1}-Q)\le dim(Q)-t,$$
\item if $(E,Q)$ is a semistable GPS with $E$ locally free and
$$rank(E_{x_1}\to Q)=a,$$
then $(E,Q)$ is $s$-equivalent to a semistable $(E',Q')$ such that 
$$dim(Tor(E')_{x_2})=dim(Q)-a.$$\endroster
The roles of $x_1$, $x_2$ in the above statements can be reversed.\endproclaim

\demo{Proof} We prove (1) at first. For given $(E',Q')\in\Cal C_{\mu}$ with
$rank(E')>0,$ there is an exact sequence
$$0\to (E'_1, Q'_1)\to (E',Q')\to (E'_2,Q'_2)\to 0$$
such that $(E'_2,Q'_2)$ is stable and $\mu_G[(E'_2,Q'_2)]=\mu$ if
$rank(E'_2)>0.$ It is clear that
$$gr(E',Q')=gr(E'_1,Q'_1)\oplus (E'_2,Q'_2).$$

When $rank(E'_2)>0,$  the $E'_2$ has to be locally free and $E'_1$ has the same
torsion with $E'$. Thus if $rank(E'_1)>0$, there is (by using induction for the rank) a $(E_1,Q_1)\in \Cal C_{\mu}$ with $E_1$ locally free and
$$rank(E_{1\,x_1}\to Q_1)\le dim(Q_1)-t$$
such that $gr(E_1,Q_1)=gr(E'_1,Q'_1).$  One can check that
$$(E,Q):=(E_1\oplus E'_2,Q_1\oplus Q'_2)\in\Cal C_{\mu}$$
is $s$-equivalent to $(E',Q')$ and
$$rank(E_{x_1}\to Q)\le dim(Q)-t.$$
If $rank(E'_1)=0$, then $gr(E',Q')=(E'_2,Q'_2)\oplus gr(Tor(E'),Tor(E')).$
Thus $(E',Q')$ satisfies (up to a $s$-equivalence) the exact sequence
$$0\to(\wt E',\wt Q')\to(E',Q')\to(\,_{x_2}\Bbb C,\Bbb C)\to 0,$$
where $(\wt E',\wt Q')\in\Cal C_{\mu}$ has torsion of dimension $t-1$ at $x_2.$
This is the typical case we treated in Lemma 2.5 of [Su], and we will indicate
later how to get our stronger statement by the construction of [Su].

When $rank(E'_2)=0$ and $dim(Tor(E'_1)_{x_2})<t$, then $(E'_2,Q'_2)$ has to be
$(\,_{x_2}\Bbb C,\Bbb C),$ which is again the above typical case we will treat.
If $dim(Tor(E'_1)_{x_2})=t$, by repeating the above procedures for $(E'_1,Q'_1),$
we will reduce the proof, after finite steps, to the above cases again since
$dim(Q'_1)$ decreases strictly. All in all, we are reduced to treating the
typical case: 
$$0\to(\wt E',\wt Q')\to(E',Q')\to(\,_{x_2}\Bbb C,\Bbb C)\to 0,$$
where $(\wt E',\wt Q')\in\Cal C_{\mu}$ and $dim(Tor(\wt E')_{x_2})=t-1.$

By using the induction for $t$, there exists
a $(\wt E,\wt Q)\in\Cal C_{\mu}$ with $\wt E$ locally free such
that
$gr(\wt E,\wt Q)=gr(\wt E',\wt Q')$ and
$$rank(\wt q_1:\wt E_{x_1}\to\wt Q)\le dim(\wt Q)-(t-1),$$
where $\tilde q_1$, $\tilde q_2$ are the induced maps by $\tilde q:\wt E_{x_1}\oplus\wt E_{x_2}\to\wt Q.$  
Since $(\,_{x_2}\Bbb C,\Bbb C)$ is stable,
we have $$gr(E',Q')=gr(\wt E,\wt Q)\oplus(\,_{x_2}\Bbb C,\Bbb C).$$
Let $K_2=ker(\tilde q_1:\wt
E_{x_2}\to\wt Q)$, choosing a Hecke modification
$h:\wt E\to E$ at $x_2$
(see Remark 1.4 of [NS]) such that
$\wt K_2:=ker(h_{x_2})\subset K_2$ and
$dim(\wt K_2)=1$, we get the
extension
$$0@>>>\wt
E@>h>>E@>\gamma>>\,_{x_2}\Bbb C@>>>0.$$
Let $Q=\wt Q\oplus\Bbb C$ and
$E_{x_2}=h_{x_2}(\wt E_{x_2})\oplus V_1$
for a subspace $V_1$, we define a
morphism $f:E_{x_1}\oplus E_{x_2}\to Q$
such that $E_{x_1}\to Q$ to be
$$E_{x_1}@>h^{-1}_{x_1}>>\wt E_{x_1}@>\tilde q_1>>\wt Q\hookrightarrow
Q$$
and $E_{x_2}\to Q$ to be
$$E_{x_2}=h_{x_2}(\wt E_{x_2})\oplus V_1@>({\bar h}^{-1}_{x_2},\gamma_{x_2})>>
\frac{\wt E_{x_2}}{\wt K_2}\oplus\Bbb
C@>(\tilde q_2,id)>>\wt Q\oplus\Bbb
C=Q$$
where $\bar h_{x_2}:\wt
E_{x_2}/\wt K_2\cong h_{x_2}(\wt E_{x_2})$ and
$\tilde q_2:\wt E_{x_2}/\wt
K_2\to\wt Q$ (note that $\wt K_2\subset K_2$).
Thus the following diagram
is commutative
$$\CD
@.\wt E_{x_1}\oplus\wt
E_{x_2}@>(h_{x_1},h_{x_2})>>E_{x_1}\oplus E_{x_2}
@>(0,\gamma_{x_2})>>\Bbb
C@>>>0\\
@.      @VqVV       @VfVV               @|           @.\\
0
@>>>\wt Q @>>> \wt Q\oplus\Bbb C  @>>> \Bbb C @>>>  0
\endCD$$
One checks
that $f$ is surjective by this diagram, and thus
$$0\to(\wt E,\wt
Q)\to(E,Q)\to(\,_{x_2}\Bbb C,\Bbb C)\to 0$$
It is easy to see that $(E,Q)\in\Cal
C_{\mu}$ is $s$-equivalent to
$(E',Q')$ and
$$rank(E_{x_1}\to Q)=rank(\wt E_{x_1}\to\wt Q)\le dim(Q)-t.$$

To prove (2), let $q:E_{x_1}\oplus E_{x_2}\to Q$ and 
$Q=q_1(E_{x_1})\oplus\Bbb C^{dim(Q)-a}$. Take the projection
$Q@>p>>\Bbb C^{dim(Q)-a}$ and define
$$\wt E:=ker(\gamma:E\to E_{x_2}@>q_2>>Q@>p>>\,_{x_2}\Bbb C^{dim(Q)-a}),$$ 
we get a semistable $(\wt E,\wt Q)\in \Cal C_{\mu}$ ($\wt Q$ being the kernel
of $p$) such that
$$0\to (\wt E,\wt Q)\to (E,Q)\to
(\,_{x_2}\Bbb C^{dim(Q)-a},\,\Bbb C^{dim(Q)-a})\to 0$$
is an exact sequence in $\Cal C_{\mu}$. Thus $(E,Q)$ is $s$-equivalent to
$$(E',Q'):=(\wt E\oplus\,_{x_2}\Bbb C^{dim(Q)-a},\wt Q\oplus\Bbb C^{dim(Q)-a})$$
by the following Lemma 2.4.
\enddemo

\proclaim{Lemma 2.4} Given a $(E,Q)\in\Cal C_{\mu}$, if there is
an exact sequence
$$0\to (E_1,Q_1)\to(E,Q)\to(E_2,Q_2)\to 0$$
such that $(E_2,Q_2)\in\Cal C_{\mu}$, then
$$gr(E,Q)=gr(E_1,Q_1)\oplus gr(E_2,Q_2).$$
In particular, $(E,Q)$ is $s$-equivalent to $(E_1\oplus E_2,Q_1\oplus Q_2).$
\endproclaim

\demo{Proof} Since $(E_2,Q_2)\in\Cal C_{\mu}$, there exists an exact sequence
$$0\to(E'_2,Q'_2)\to(E_2,Q_2)\to(E''_2,Q''_2)\to 0$$
such that $(E''_2,Q''_2)\in\Cal C_{\mu}$ is stable. Thus
$$gr(E_2,Q_2)=gr(E'_2,Q'_2)\oplus(E''_2,Q''_2).$$
On the other hand, if we define $(\wt E,\wt Q)$ by exact sequence
$$0\to(\wt E,\wt Q)\to(E,Q)@>g>>(E''_2,Q''_2)\to 0,$$
where $g:(E,Q)\to(E_2,Q_2)\to(E''_2,Q''_2),$ then we have an exact sequence
$$0\to (E_1,Q_1)\to(\wt E,\wt Q)\to(E'_2,Q'_2)\to 0,$$
and $(E'_2,Q'_2)\in\Cal C_{\mu}$. By using the induction for the 
$rank(E_2)$ and $h^0(E_2)$ when $rank(E_2)=0,$ we have
$$gr(\wt E,\wt Q)=gr(E_1,Q_1)\oplus gr(E'_2,Q'_2).$$
Now the lemma is clear.\enddemo

\heading \S3 The factorization theorem \endheading

Recall that $\pi:\wt X_0\to X_0$ is the normalisation of $X_0$ and
$\pi^{-1}(x_0)=\{x_1,x_2\}$. Given a GPS $(E, E_{x_1}\oplus E_{x_2}@>q>>Q)$
on $\wt X_0$, we define a coherent sheaf $\phi(E,Q):=F$ by the exact
sequence $$0\to F\to\pi_*(E)\to\, _{x_0}Q\to 0,$$
where we use $_xW$ to denote the skyscraper sheaf supported at $\{x\}$ with fibre $W$, and the morphism $\pi_*(E)\to\, _{x_0}Q$ is defined as follows
$$\pi_*(E)\to\pi_*(E)|_{\{x_0\}}=\,_{x_0}(E_{x_1}\oplus E_{x_2})@>q>>\, _{x_0}Q.$$
It is clear that $F$ is torsion free of rank $(r_1,r_2)$ if and only if $(E,Q)$
is a GPS of rank $(r_1,r_2)$ and satisfying 
$$(\text{Tor}E)_{x_1}\oplus(\text{Tor}E)_{x_2}\overset q
\to\hookrightarrow Q.\tag T$$
In particular, the GPS in $\SH$ give in this way torsion free sheaves of 
rank $r$ with the natural parabolic structures at points of $I$.

\proclaim{Lemma 3.1} Let $(E,Q)$ satisfy
condition (T), and $F=\phi(E,Q)$
the associated torsion free sheaf on $X_0$.
We have
\roster\item If $E$ is a vector bundle and the maps $E_{x_i}\to
Q$
are isomorphisms, then $F$ is a vector bundle.
\item If $F$ is a vector
bundle on $X_0$, then there is an unique $(E,Q)$
such that $\phi(E,Q)=F.$ In
fact, $E=\pi^*F.$
\item If $F$ is a torsion free sheaf, then there is a
$(E,Q)$, with
$E$ a vector bundle on $\wt X_0,$ such that $\phi(E,Q)=F$ and
$E_{x_2}\to Q$
is an isomorphism. The rank of the map $E_{x_1}\to Q$ is $a$
iff 
$F\otimes\hat\SO_{x_0}\cong \hat\SO_{x_0}^{\oplus a}\oplus
m_{x_0}^{\oplus
(r-a)}.$ 
The roles of $x_1$ and $x_2$ can be
reversed.
\item Every torsion free rank $r$ sheaf $F$ on $X_0$ comes from a
$(E,Q)$
such that $E$ is a vector
bundle.\endroster\endproclaim

\demo{Proof} Similar with Lemma 4.6 of [NR] and Lemma 2.1 of [Su].
\enddemo

\proclaim{Lemma 3.2} Let $F=\phi(E,Q)$,
then $F$ is semistable if and
only if $(E,Q)$ is semistable. Moreover, one
has
\roster\item If $(E,Q)$ is stable, then $F$ is stable.
\item If $F$ is
a stable vector bundle, then $(E,Q)$ is
stable.\endroster
\endproclaim
\demo{Proof} For any subsheaf $E'\subset E$ such that $E/E'$ is torsion free 
outside $\{x_1,x_2\}$, the induced GPS $(E',Q^{E'})$ defines a subsheaf
$F'\subset F$ by
$$ 0\to F'\to\pi_*(E')\to\, _{x_0}Q^{E'}\to 0.$$
It is clear that $par\chi_m(F')=par\chi_m(E')-dim(Q^{E'})$, thus $F$ semistable
implies $(E,Q)$ semistable. Note that $E$ may have torsion and thus 
$(E,Q)$ may not
be stable even if $F$ is stable (for instance, taking $E'$ to be the torsion 
subsheaf). In fact, $(E,Q)$ is stable if and only if $F$ is a stable vector bundle.  

Next we prove that if $(E,Q)$ is stable (semistable), then $F$ is stable
(semistable). For any subsheaf $F'\subset F$ such that $F/F'$ is torsion free,
we have canonical morphism $\pi^*F'\to \pi^*F\to\pi^*\pi_*E\to E$. Let
$E'$ be the image of $\pi^*F'$, one has the diagram 
$$\CD
@.      0                @.        0  @.     0      \\
@.    @VVV                    @VVV    @VVV     @. \\
0 @>>> F'          @>>> \pi_*E'@>>>\,_{x_0}Q^{E'}@>>>0  \\
@.      @VVV             @VVV         @VVV      @.\\
0 @>>>F  @>>>\pi_*E @>>>\,_{x_0}Q@>>>0       \\
@. @VVV             @VVV          @VVV       @.\\
0 @>>> F/F'@>>>\pi_*(E/E')@>>>\,_{x_0}(Q/Q^{E'})@>>>0   \\
@.      @VVV  @VVV         @VVV         @.\\
@.      0    @.           0   @.        0 @.  @.
\endCD$$
which implies $E/E'$ torsion free outside $\{x_1,x_2\}$ (since $F/F'$
torsion free) and 
$$par\chi_m(F')=par\chi_m(E')-dim(Q^{E'}),\quad par\chi_m(F)=par\chi_m(E)-dim(Q).$$
Thus, note that $rk(E')=rk(F')$ and $rk(E)=rank(F)$, one proves the lemma.   
\enddemo

\proclaim{Lemma 3.3} Let $(E,Q)$ be a semistable GPS with $E$ locally free
and $F=\phi(E,Q)$ the associated torsion free sheaf, if there exists an exact
sequence $$0\to F_1\to F\to F_2\to 0$$
with $F_2$ semistable and $par\mu_m(F_2)=par\mu_m(F)$. Then
$(E,Q)$ is $s$-equivalent to a semistable $(E',Q')$ such that $E'$ has torsion of dimension $$dim(Q)-\bold a(F_1)-\bold a(F_2).$$\endproclaim

\demo{Proof} For any torsion free sheaf $F$, we have a canonical exact
sequence $$0\to F\to \pi_*\wt E\to \wt Q\to 0$$
where $\wt E=\pi^*F/Tor(\pi^*F)$ and $dim(\wt Q)=\bold a(F)$. If $F=\phi(E,Q)$
with $E$ locally free, then we have a commutative diagram
$$\CD
@.        @.           0   @.        0  @.           \\
@.
@.                     @VVV          @VVV         @. \\
0  @>>> F         @>>>\pi_*\wt E @>>>\,_{x_0}\wt Q    @>>> 0  \\
@.      @|             @VVV         @VVV        @.  \\
0 @>>> F @>>> \pi_*E  @>>>\,_{x_0}Q  @>>>0   \\
@.    @.             @VVV          @VVV        @. \\
@.       @.       \pi_*\tau @>>>\,_{x_0}Q'      @.   \\
@.      @.              @VVV         @VVV        @.\\
@.         @.           0   @.        0   @.
\endCD$$
where $\tau=E/\wt E$ and $Q'=Q/\wt Q$, the map $\pi_*\tau\to\,_{x_0}Q'$ 
is defined
such that the diagram is commutative, which has to be an isomorphism.
This gives an exact sequence 
$$0\to (\wt E,\wt Q)\to (E,Q)\to (\tau,Q')\to 0$$
in $\Cal C_{\mu}$, thus $(E,Q)$ is $s$-equivalent to
$(\wt E\oplus\tau,\wt Q\oplus Q').$  On the other hand, we consider the
following commutative diagram
$$\CD
@.      0  @.           0   @.        0  @.           \\
@.
@VVV            @VVV          @VVV         @. \\
0  @>>> F_1         @>>> F
@>>> F_2     @>>> 0  \\
@.      @VVV             @VVV         @VVV
@.  \\
0 @>>>\pi_* E_1  @>>>   \pi_*\wt E@>>>\pi_*E_2      @>>>0   \\
@.
@VVV             @VVV          @VVV        @.
\\
0@>>>\,_{x_0} Q_1@>>>\,_{x_0}\wt Q@>>>\,_{x_0}Q_2@>>> 0    \\
@.      @VVV
@VVV         @VVV        @.\\
@.      0    @.           0   @.        0
@.
\endCD\tag3.1$$ 
where $E_1=\pi^*F_1/Tor(\pi^*F_1)$, $dim(Q_1)=\bold a(F_1)$, the first two
vertical sequences are the canonical exact sequences determined by 
$F_1$ and $F$, and $E_2=\wt E/E_1$, $Q_2=\wt Q/Q_1$, the third vertical sequence
is defined by demanding the diagram commutative, which has to be exact. It is
easy to see that $\mu_G[(E_2,Q_2)]=\mu_G[(\wt E,\wt Q)]$ and $(E_2,Q_2)$
is semistable (since $F_2$ is so). Thus 
$$gr(\wt E,\wt Q)=gr(E_1,Q_1)\oplus gr(E_2,Q_2),$$
which implies that $(E,Q)$ is $s$-equivalent to $$(E',Q'):=(E_1\oplus E_2\oplus\tau,Q_1\oplus Q_2\oplus Q').$$
One checks that $dim(Tor(E_2))=\bold a(F)-\bold a(F_1)-\bold a(F_2)$ by
restricting the diagram (3.1) to point $x_0$ and counting the dimension
of fibres (the first two vertical sequences remaining exact). Therefore
$$dim(Tor(E'))=dim(\tau)+dim(Tor(E_2))=dim(Q)-\bold a(F_1)-\bold a(F_2),$$
we have proved the lemma.
\enddemo

Consider the family $\rho^*\SE=(\rho^*\SE^1,\rho^*\SE^2)$ of GPS over $\wt\SR^{ss}$ with the universal
quotient $\rho^*(\SE^1_{x_1}\oplus\SE^2_{x_2})\to\SQ,$ using the finite morphism
$$\pi\times I_{\wt\SR^{ss}}:\wt X_0\times\wt\SR^{ss}\to X_0\times\wt\SR^{ss},$$
we can define a family $\SF_{\wt\SR^{ss}}$ of semistable sheaves (Lemma 3.2) on 
$X_0$ by the exact sequence
$$0\to\SF_{\wt\SR^{ss}}\to(\pi\times I_{\wt\SR^{ss}})_*(\rho^*\SE)\to\,_{x_0}\SQ
\to 0\tag3.2.$$
Since $\rho^*\SE$ is flat over $\wt\SR^{ss}$ and $\SQ$ locally free on
$\wt\SR^{ss}$, $\SF_{\wt\SR^{ss}}$ is a flat family over $\wt\SR^{ss}$. Thus
we have a morphsim
$$\phi_{\wt\SR^{ss}}:\wt\SR^{ss}\to\SR^{ss}\to\SU_{X_0}$$
such that $\phi_{\wt\SR^{ss}}^*\Theta_{\SU_{X_0}}=\Theta_{\SF_{\wt\SR^{ss}}}$
by Theorem 1.2 in the $\S1$. 

\proclaim{Lemma 3.4} The morphism $\phi_{\wt\SR^{ss}}$ induces a
morphism $$\phi_{\SP_{\chi_1,\chi_2}}:\SP_{\chi_1,\chi_2}\to \SU_{X_0}$$
such that $\phi_{\SP_{\chi_1,\chi_2}}^*\Theta_{\SU_{X_0}}=\Theta_{\SP_{\chi_1,\chi_2}}.$
\endproclaim

\demo{Proof} The proof is clear, we just remark that one can compute 
 $\Theta_{\SF_{\wt\SR^{ss}}}=\Theta_{\wt\SR^{ss}}$ by the exact sequence (3.2)
defining the sheaf $\SF_{\wt\SR^{ss}}.$\enddemo

Let $\SU_{\chi_1,\chi_2}$ be the image of $\SP_{\chi_1,\chi_2}$ under the
morphism $\phi_{\SP_{\chi_1,\chi_2}}$, then $\SU_{\chi_1,\chi_2}$ is an 
irreducible component of $\SU_{X_0}$ and $\phi_{\SP_{\chi_1,\chi_2}}$ is a finite morphism since it pulls back an ample line bundle to an ample line
bundle. We will see that   
$$\phi_{\SP_{\chi_1,\chi_2}}:\SP_{\chi_1,\chi_2}\ssm\{\SD_1,\SD_2\} \to 
\SU_{\chi_1,\chi_2}\ssm\SW_{r-1}$$ is an isomorphism. Thus
$\phi_{\SP_{\chi_1,\chi_2}}$ is the normalisation of $\SU_{\chi_1,\chi_2}.$
We have clearly the morphism
$$\phi:=\coprod_{\chi_1+\chi_2=\chi+r}\phi_{\SP_{\chi_1,\chi_2}}:
\SP\to\SU_{X_0},$$
which is the normalisation of $\SU_{X_0}.$ We copy Proposition 2.1 from
[Su]. 

\proclaim{Proposition 3.1} With the above notation
and denoting
$\SD_1(r-1)$,
$\SD_2(r-1)$, $\SW_{r-1}$ by $\SD_1$, $\SD_2$
and $\SW$, we have
\roster\item $\phi: \SP\to\SU_{X_0}$ is finite and
surjective, and
$\phi(\SD_1(a))=\phi(\SD_2(a))=\SW_a$,
\item $\phi(\SP\ssm
\{\SD_1\cup\SD_2\})=\SU_{X_0}\ssm\SW$, and induces an
isomorphism on
$\SP\ssm\{\SD_1\cup\SD_2\},$
\item $\phi |_{\SD_1(a)}:\SD_1(a)\to\SW_a$ is
finite and surjective,
\item
$\phi(\SD_1(a)\ssm\{\SD_1(a)\cap\SD_2\cup\SD_1(a-1)\})
=\SW_a\ssm\SW_{a-1}$, and
induces an isomorphism on
$\SD_1(a)\ssm\{
\SD_1(a)\cap\SD_2\cup\SD_1(a-1)\},$
\item
$\phi:\SP\to\SU_{X_0}$ is the normalisation of $\SU_{X_0}$,
\item $\phi|_{\SD_1(a)}:\SD_1(a)\to\SW_a$ is the normalisation
of $\SW_a$,
\item $\phi(\SD_1(a)\cap\SD_2)=\SW_{a-1}$, and $\SW_{a-1}$ is the non-normal
locus of $\SW_a$.
\endroster\endproclaim

\demo{Proof} In proving (4), we used Lemma 2.6 of [Su] to show that
 $\phi$ induces a morphism
$$\phi:\SD_1(a)\ssm\{\SD_1(a)\cap\SD_2\cup\SD_1(a-1)\}
@>>>\SW_a\ssm\SW_{a-1}.$$
But Lemma 2.6 in [Su] is not correct, we have to prove it without using
the lemma (also to fix the gap in [Su]). We will use $[\quad]$ to denote
the $s$-equivalent class of objects we are considering. For any
$[(E,Q)]\in\SD_1(a)\ssm\{\SD_1(a)\cap\SD_2\cup\SD_1(a-1)\}$,
we
can assume that $E$ is a vector bundle by Proposition 2.5, and $E_{x_2}\to
Q$
is an isomorphism since $[(E,Q)]\notin \SD_2$. Thus
$\phi(E,Q)=F
\in \hat\SW_a\ssm\hat\SW_{a-1}$ by Lemma 3.1 (3), we need to show
$[F]\notin\SW_{a-1}$. If it is not so, then $F$ is $s$-equivalent to
a semistable torsion free sheaf $F'\in\hat\SW(a-1)$ and has an exact
sequence $$0\to F_1\to F\to F_2\to 0$$
with $par\mu_m(F_2)=par\mu_m(F)$ and $F_2$ stable. Thus 
$gr(F')=gr(F_1)\oplus F_2$ and (by Lemma 1.3)
$$a-1\ge\bold a(F')\ge \bold a(F_1)+\bold a(F_2).$$
On the other hand, by Lemma 3.3, $(E,Q)$ is $s$-equivalent to a semistable
$(E',Q')$ with $dim(Tor(E'))=r-\bold a(F_1)-\bold a(F_2).$ By Proposition 2.5 (1), $E'$ has no torsion at $x_1$ since $[E',Q')]=[(E,Q)]\notin\SD_2$. Hence,
by Proposition 2.5 (1) again,
$(E',Q')$ is $s$-equivalent to a $(\wt E,\wt Q)$ with $\wt E$ locally free
and $$rank(\wt E_{x_1}\to\wt Q)\le\bold a(F_1)+\bold a(F_2)\le a-1,$$
we get the contradiction $[(E,Q)]=[(\wt E,\wt Q)]\in \SD_1(a-1).$
Thus $\phi$ induces a morphism
$$\phi:\SD_1(a)\ssm\{\SD_1(a)\cap\SD_2\cup\SD_1(a-1)\}
@>>>\SW_a\ssm\SW_{a-1}$$

The argument in [Su] for other statements goes through, only (7) is in doubt. This can
be seen as follows, the fact 
$\phi(\SD_1(a)\cap\SD_2)=\SW_{a-1}$ follows the local compuation (see Proposition 3.9 of [B3]), and the non-normal locus of $\SW_a$ is contained in
$\SW_{a-1}$ by the above (4). If $\SW_{a-1}$ is not empty and not equal to
the non-normal locus, there exists a non-empty irreducible component 
$\SW_{a-1}^{\chi_1,\chi_2}$ of $\SW_{a-1}$ such that $\phi |_{\SD_1(a)}$ is
an isomorphism at the generic point of $\SW_{a-1}^{\chi_1,\chi_2}.$ It
is impossible since the fibre has at least two points (one is in 
$\SD_1(a-1)\ssm\SD_2$ by (3) of Lemma 3.1, another is in $\SD_1(a)\cap\SD_2$).
\enddemo

Let $I_Z$  denote the ideal sheaf of closed subscheme $Z$ in a scheme $X$. 
When $Z$ is of codimension one (not necessarily a Cartier divisor), 
we set $\SO_X(-Z):=I_Z$. If $\SL$ is a line bundle on $X$ and $Y$ is a closed
subscheme of $X$, we denote $\SL\otimes I_Z$ and the restriction $I_Z\otimes
\SO_Y$ of $I_Z$ on $Y$ by $\SL(-Z)$ and $\SO_Y(-Z)$. We have the straightforward generalisations of [Su, Lemma 4.3 and Proposition 4.1], whose proof we omit.

\proclaim{Lemma 3.5} Suppose given a seminormal variety $V$, with
normalization $\sigma:\wt V\to V$. Let the non-normal locus be $W$,
endowed with its reduced structure. Let $\wt W$ be set-theoretic inverse
image of $W$ in $\wt V$, endowed with its reduced structure. Let
$N$ be a line bundle on $V$, and let $\wt N$ be its pull-back to
$\wt V$ ($\wt N=\sigma^*N$). Suppose $H^0(\wt V,\wt N)\to H^0(\wt W,\wt N)$
is surjective. Then
\roster
\item 
There is an exact sequence
$$0\to H^0(\wt V,\wt N\otimes I_{\wt W})\to H^0(V,N)\to H^0(W,N)\to 0.$$
\item If $H^1(W,N)\to H^1(\wt W,\wt N)$ is injective, so is
$H^1(V,N)\to H^1(\wt V,\wt N).$
\endroster\endproclaim

\proclaim{Lemma 3.6} The following maps are surjective for any
$1\le a\le r$
\roster \item $H^0(\SD_1(a),\Theta_{\SP})\to H^0(\SD_1(a)\cap\SD_2\cup\SD_1(a-1),\Theta_{\SP})$.
\item $H^0(\SD_1(a),\Theta_{\SP})\to H^0(\SD_1(a)\cap\SD_2,\Theta_{\SP}).$\endroster\endproclaim

The above Lemma 3.6 tells us that the assumption (surjectivity) in Lemma 3.5
is satisfied for the situation: $V=\SW_a$, $\wt V=\SD_1(a),$
$\sigma=\phi|_{\SD_1(a)}$ and $N=\Theta_{\SU_X}|_{\SW_a}.$ Thus we can use
Lemma 3.5 to prove that

\proclaim{Proposition 3.2} We have a (noncanonical) isomorphism 
$$H^0(\SU_{X_0},\Theta_{\SU_{X_0}})\cong H^0(\SP,\Theta_{\SP}(-\SD_2)).$$
\endproclaim

\demo{Proof} Similar with the proof of Proposition 4.3 of [Su].\enddemo

\proclaim{Proposition 3.3} Let $\wt\SR_F\subset\SH$ be the open set consisting
of $(E,Q)$ with $E$ locally free. Then
$$H^0(\wt\SR^{ss},\Theta_{\wt\SR^{ss}})^G=H^0(\SH,\Theta_{\SH})^G=
H^0(\wt\SR_F,\Theta_{\wt\SR_F})^G,$$
where $G=(GL(V_1)\times GL(V_2))\cap SL(V_1\oplus V_2).$\endproclaim

\demo{Proof} The first equality follows the following Lemma 3.7, the
second equality follows the following Lemma 3.8 by taking
$V=\wt\SR^{ss},$  $U=\SP_{\chi_1,\chi_2},$ $V'=\wt\SR^{ss}\cap\wt\SR_F$
and $U^{\prime\prime}=\SP_{\chi_1,\chi_2}\ssm\{\SD_1,\SD_2\}$ (one need
here Proposition 1.4 to show that $U^{\prime\prime}$ is nonempty).\enddemo

\proclaim{Lemma 3.7} Let $V$ be a projective scheme on which a reductive
group $G$ acts, $\wt\SL$ an ample line bundle linearising the $G$-action,
and $V^{ss}$ the open subscheme of semistable points. Let $V'$ be a
$G$-invariant closed subscheme of $V^{ss}$, $\bar V'$ its schematic closure
in $V$. Then
\roster\item $\bar V^{'ss}=V'$, and $V'\diagup\diagup G$ is a closed subscheme
of $V^{ss}\diagup\diagup G.$
\item $H^0(V^{ss},\wt\SL)^{inv}=H^0(W,\wt\SL)^{inv}$, where $W$ is an open 
$G$-invariant irreducible normal subscheme of $V$ containing $V^{ss}$ and
$(\quad)^{inv}$ denotes the invariant subspace for an action of $G$.
\endroster\endproclaim

\demo{Proof}  See Lemma 4.14 and Lemma 4.15 of [NR].\enddemo

\proclaim{Lemma 3.8} Let $V$ be a normal variety with a $G$-action, where $G$
is a reductive algebraic group. Suppose a good quotient $\pi:V\to U$ exists.
Let $\wt\SL$ be a $G$-line bundle on $V$, and suppose it descends as a line
bundle $\SL$ on $U$. Let $V^{\prime\prime}\subset V'\subset V$ be open $G$-invariant 
subvarieties of $V$, such that $V'$ maps onto $U$ and 
$V^{\prime\prime}=\pi^{-1}(U^{\prime\prime})$ for
some nonempty open subset $U^{\prime\prime}$ of $U$. Then any invariant section of
$\wt\SL$ on $V'$ extends to $V$.\endproclaim

\demo{Proof}  See Lemma 4.16 of [NR].\enddemo

\proclaim{Proposition 3.4} Let $G_1$ and $G_2$ be reductive algebraic 
groups acting on the normal projective schemes $\bar V_1$, $\bar V_2$ 
with ample linearizing $L_1$ and $L_2$.
Suppose that $L_1$ and $L_2$ descend to $\Theta_1$ and $\Theta_2$. Then,
for any open sets $V_1\supset \bar V_1^{ss}$ and $V_2\supset\bar V_2^{ss}$,
$$H^0(V_1\times V_2, L_1\otimes L_2)^{G_1\times G_2}=H^0(V_1,L_1)^{G_1}\otimes
H^0(V_2,L_2)^{G_2}.$$\endproclaim

\demo{Proof} Using Lemma 3.7 and next Lemma 3.9, we have
$$\aligned &H^0(V_1\times V_2, L_1\otimes L_2)^{G_1\times G_2}\\
&=\left\{H^0(V_1\times V_2, L_1\otimes L_2)^{G_1\times\{id\}}\right\}^{\{id\}\times G_2}\\
&=H^0(\bar V_1^{ss}//G_1\times V_2, \Theta_1\otimes L_2)^{\{id\}\times G_2}\\
&=H^0(\bar V_1^{ss}//G_1\times\bar V_2^{ss}//G_2, \Theta_1\otimes\Theta_2)\\
&=H^0(\bar V_1^{ss}//G_1,\Theta_1)\otimes H^0(\bar V_2^{ss}//G_2,\Theta_2)\\
&=H^0(V_1,L_1)^{G_1}\otimes H^0(V_2,L_2)^{G_2}.\endaligned$$\enddemo

\proclaim{Lemma 3.9} Suppose $V\to V//G$ is a good quotient and $T$ is any
 variety with trivial $G$-action. Then
$V\times T\to V//G\times T$ is a good quotient.\endproclaim

\proclaim{Notation 3.1}For $\mu=(\mu_1,\cdots,\mu_r)$  with $0\le\mu_r\le\cdots\le\mu_1\le k-1,$
let $$\{d_i=\mu_{r_i}-\mu_{r_i+1}\}_{1\le i\le l}$$ be the subset of nonzero
integers in $\{\mu_i-\mu_{i+1}\}_{i=1,\cdots,r-1}.$ Then we define that
$$ r_i(x_1)=r_i,\quad d_i(x_1)=d_i,\quad l_{x_1}=l,\quad \alpha_{x_1}=\mu_r$$
$$ r_i(x_2)=r-r_{l-i+1},\quad d(x_2)=d_{l-i+1},\quad l_{x_2}=l,\quad\alpha_{x_2}=k-\mu_1$$
and for $j=1,2$, we set
$$\aligned
\vec a(x_j)&=\left(\mu_r,\mu_r+d_1(x_j),\cdots,\mu_r+
\sum^{l_{x_j}-1}_{i=1}d_i(x_j),\mu_r+\sum^{l_{x_j}}_{i=1}d_i(x_j)\right)\\ 
\vec n(x_j)&=(r_1(x_j),r_2(x_j)-r_1(x_j),
\cdots,r_{l_{x_j}}(x_j)-r_{l_{x_j}-1}(x_j)).\endaligned$$
We also define that
$$\aligned&\chi_1^{\mu}=\frac{1}{k}\left(\sum_{x\in I_1}\sum^{l_x}_{i=1}
d_i(x)r_i(x)+r\sum_{x\in I_1}\alpha_x+r\ell_1\right)+\frac{1}{k}\sum^r_{i=1}\mu_i\\
&\chi_2^{\mu}=\frac{1}{k}\left(\sum_{x\in I_2}\sum^{l_x}_{i=1}
d_i(x)r_i(x)+r\sum_{x\in I_2}\alpha_x+r\ell_2\right)+r-\frac{1}{k}\sum^r_{i=1}\mu_i.\endaligned$$
\endproclaim

One can check that the numbers defined in Notation 3.1 satisfy ($j=1,2$)
$$ \sum_{x\in
I_j\cup\{x_j\}}\sum^{l_x}_{i=1}d_i(x)r_i(x)+r\sum_{x\in
I_j\cup\{x_j\}}
\alpha_x+r\ell_j=k\chi_j^{\mu}.\tag3.3$$

\proclaim{Notation 3.2} For the numbers defined in Notation 3.1, let, for $j=1,2$, 
$$\SU^{\mu}_{X_j}:=\SU_{X_j}(r,\chi_j^{\mu}, I_j\cup\{x_j\},
\{\vec n(x),\vec a(x)\}_{x\in I_j\cup\{x_j\}},k)$$
be the moduli space of $s$-equivalence classes of semistable
parabolic bundles $E$ of rank $r$ on $X_j$ and $\chi(E)=\chi_j^{\mu}$,
together with parabolic structures of type 
$\{\vec n(x)\}_{x\in I\cup\{x_j\}}$ and weights 
$\{\vec a(x)\}_{x\in I\cup\{x_j\}}$ at points
$\{x\}_{x\in I\cup\{x_j\}}$. We define $\SU^{\mu}_{X_j}$ to be empty 
if $\chi_j^{\mu}$ is not an integer. Let
$$\Theta_{\SU^{\mu}_{X_j}}:=\Theta(k,\ell_j,\{\vec n(x),\vec a(x),\alpha_x\}_{x\in I_j\cup\{x_j\}}, I_j\cup\{x_j\})$$
be the theta line bundle.\endproclaim

\proclaim{Theorem 3.1} There exists a (noncanonical) isomorphism
 $$H^0(\SU_{X_1\cup X_2},\Theta_{\SU_{X_1\cup X_2}})\cong\bigoplus_{\mu}
H^0(\SU^{\mu}_{X_1},\Theta_{\SU^{\mu}_{X_1}})\otimes
H^0(\SU^{\mu}_{X_2},\Theta_{\SU^{\mu}_{X_2}})$$
where $\mu=(\mu_1,\cdots,\mu_r)$ runs through the integers $0\le\mu_r\le\cdots\le
\mu_1\le k-1.$\endproclaim

\demo{Proof} By Proposition 3.3, one can show shat
$$H^0(\SP_{\chi_1,\chi_2},\Theta_{\SP_{\chi_1,\chi_2}}(-\SD_2))=
H^0(\wt\SR_F,\Theta_{\wt\SR_F}(-\hat\SD_2))^G.$$
Note that $\SO_{\wt\SR_F}(-\hat\SD_2)=det\SE_{x_2}\otimes (det\SQ)^{-1}$ and
write $\eta_{x_2}:=(det\SE_{x_2})^{-1}\otimes det\SQ$, we have 
$$H^0(\wt\SR_F,\Theta_{\wt\SR_F}(-\hat\SD_2))^G=
H^0(\SR_{1\,F}\times\SR_{2\,F},\Theta_{\SR_{1\,F}}\otimes\Theta_{\SR_{2\,F}}
\otimes(det\SE_{x_2})^k\otimes\rho_*(\eta_{x_2}^{k-1}))^G.$$
Let 
$$\SR_j^{\mu}:=\underset{x\in I_j\cup\{x_j\}}\to{\times_{\bold{\wt Q_{1\,F}}}}
Flag_{\vec n(x)}(\Cal F^j_x)@>p^{\mu}_j>>\SR_{j,F},$$
then, by Lemma 4.6 of [Su], we have
$$\rho_*(\eta_{x_2}^{k-1})=\bigoplus_{\mu}p^{\mu}_{1*}(\SL^{\mu}_1)\otimes
p^{\mu}_{2*}(\SL^{\mu}_2)$$
where $\mu=(\mu_1,\cdots,\mu_r)$ runs through the integers 
$0\le\mu_r\cdots\le\mu_1\le k-1$ and
$$\SL_1^{\mu}=(det\,\SE^1_{x_1})^{\mu_r}\otimes
\bigotimes^{l_{x_1}}_{i=1}(det\SQ_{x_1,i})^{d_i(x_1)},$$ $$\SL_2^{\mu}=(det\,\SE^2_{x_2})^{-\mu_1}\otimes
\bigotimes^{l_{x_2}}_{i=1}(det\SQ_{x_2,i})^{d_i(x_2)}$$
are line bundles on $\SR_1^{\mu}\times\SR_2^{\mu}.$ By the definition
$$\Theta_{\SR_j^{\mu}}:=(detR\pi_{\SR^{\mu}_j}\SE^j)^k\otimes\bigotimes_{x\in I_j\cup\{x_j\}}\left\{(det
\SE^j_x)^{\alpha_x}\otimes\bigotimes^{l_x}_{i=1}
(det\SQ_{x,i})^{d_i(x)}\right\}\otimes(det\SE^j_{y_1})^{\ell_j},$$
one sees easily that 
$$\Theta_{\SR_1^{\mu}}=p^{\mu\,*}_1(\Theta_{\SR_{1\,F}})\otimes\SL_1^{\mu},$$
$$\Theta_{\SR_2^{\mu}}=p^{\mu\,*}_2(\Theta_{\SR_{2\,F}}\otimes (det\SE_{x_2})^k)\otimes\SL_2^{\mu}.$$
Thus we have (for any $\chi_1,\chi_2$) the equality
$$H^0(\SP_{\chi_1,\chi_2},\Theta_{\SP_{\chi_1,\chi_2}}(-\SD_2))=
H^0(\SR_1^{\mu}\times\SR_2^{\mu},
\Theta_{\SR_1^{\mu}}\otimes\Theta_{\SR_2^{\mu}})^G.$$
Since $\Bbb C^*\times\Bbb C^*$ acts trivially on $\SR_1^{\mu}\times\SR_2^{\mu}$,
one can see that if 
$$H^0(\SR_1^{\mu}\times\SR_2^{\mu},
\Theta_{\SR_1^{\mu}}\otimes\Theta_{\SR_2^{\mu}})^G\neq 0,$$ then the $\chi_j$
($j=1,2$) has to satisfy
$$ \sum_{x\in
I_j\cup\{x_j\}}\sum^{l_x}_{i=1}d_i(x)r_i(x)+r\sum_{x\in
I_j\cup\{x_j\}}
\alpha_x+r\ell_j=k\chi_j.$$
Therefore $\chi_j$ has to be $\chi_j^{\mu}$. In this case,  $\Bbb C^*\times\Bbb C^*$ acts trivially on the line bundle, 
$$H^0(\SR_1^{\mu}\times\SR_2^{\mu},
\Theta_{\SR_1^{\mu}}\otimes\Theta_{\SR_2^{\mu}})^G=
H^0(\SR_1^{\mu}\times\SR_2^{\mu},\Theta_{\SR_1^{\mu}}\otimes\Theta_{\SR_2^{\mu}})^{SL(V_1)\times SL(V_2)}.$$
Thus, by using Proposition 3.4, we can prove the theorem.\enddemo

We end this paper by some remarks. In Notation 1.1, we chose and fixed the
ample line bundle $\SO(1),$ the theta line bundle and the factorization
are generally depend on this choice. In some cases, although the moduli space
itself depends the choice, the theta bundle and the factorization (also
the number of irreducible components of the moduli space) are independent of
the choice. For example, when $\chi=0$, $|I|=0$, or the parabolic degree is zero, we can manage to the case: $\ell_1+\ell_2=0.$  In any case, one can see
that $\chi_1^{\mu}<n_1+r$, thus, for any choice, there are only $r$ components
of moduli space contribute to the factorization.

The choice in Notation 1.1 has quit freedom, it is in general a choice of
the partitions of $\ell_1+\ell_2$. In particular, if we are only interested
in studying moduli space, we can choose any $\SO(1)$.

\proclaim{Corollary 3.1} There is a choice of $\SO(1)$ such that the moduli
space $\SU_{X_0}$ has  $r$ irreducible components and
$$\SW_0=\emptyset.$$
In particular, when $r=2$, $\SU_{X_0}$ has two normal crossing irreducible
components.\endproclaim

\demo{Proof} One can easily choose $\SO(1)$ such that $n_1$ and $n_2$ are not
integers. Thus $n_j<\chi_j<n_j+r$ ($j=1,2$) has only $r$ possibility and each
such $\chi_j$ there is a nonempty irreducible component by Proposition 1.4.
Recall (2.1)
$$\aligned&n_1+r-dim(Q^{E_2})\le\chi(E_1)\le n_1+dim(Q^{E_1})\\
&n_2+r-dim(Q^{E_1})\le\chi(E_2)\le n_2+dim(Q^{E_2}),\endaligned$$
we see that $dim(Q^{E_j})\ge \chi_j-n_j>0$, which means that
$$\SD_1(0)= \SD_2(0)=\emptyset.$$
Thus $\SW_0=\emptyset.$ In particular, when $r=2$, the local model of moduli
space at any non-locally free sheaf is $\Bbb C[x,y]/(xy)$ by Lemma 1.4.
\enddemo

\remark{Remark 3.1} When $r=2$ and $\SO(1)$ is chosen such that $n_1$ and $n_2$ are 
not integers, $\SP$ has two disjoint irreducible components $\SP_1$ and $\SP_2$,
$\SD_j\subset\SP_j$ ($j=1,2$) is isomorphic to $\SW\subset\SU_{X_0}.$ Thus
$\SU_{X_0}$ can be obtained from $\SP_1$ and $\SP_2$ by identifying $\SD_1$
and $\SD_2$.\endremark


\Refs

\widestnumber\key{B1}
\widestnumber\key{B2}
\widestnumber\key{BR}
\widestnumber\key{CS}
\widestnumber\key{DN}
\widestnumber\key{Fa}
\widestnumber\key{La}
\widestnumber\key{MS}
\widestnumber\key{Ne}
\widestnumber\key{NR}
\widestnumber\key{NS}
\widestnumber\key{Pa}
\widestnumber\key{S1}
\widestnumber\key{S2}
\widestnumber\key{Si}
\widestnumber\key{Tr}
\widestnumber\key{EGA-I}

\ref\key Be\by A. Beauville\paper Vector bundles on curves and generalized
theta functions: recent results and open problems.\jour Math. Sci. Res. Inst.
Publ.\vol 28\pages 17--33\yr1995\endref

\ref\key B1 \by U. Bhosle\paper Generalised parabolic bundles and
applications to torsionfree sheaves on nodal curves
\pages 187--215\yr1992\vol 30 \jour Arkiv f{\"o}r matematik\endref

\ref\key B2 \by U. Bhosle\paper Generalized parabolic bundles and
applications--II\pages 403--420 \yr1996\vol 106 \jour Proc. Indian
Acad. Sci. (Math. Sci.)\endref

\ref\key B3 \by U. Bhosle\paper Vector bundles on curves with many 
components\pages 81--106\yr1999\vol 79\jour Proc. London Math. Soc.
\endref

\ref\key DW1\by G. Daskalopoulos and R. Wentworth\paper Local degeneration
of the moduli space of vector bundles and factorization of rank two theta
functions. I\jour Math. Ann.\vol 297\pages 417--466\yr1993\endref

\ref\key DW2\by G. Daskalopoulos and R. Factorization of rank two theta
functions. II: Proof of the Verlinde formula\jour Math. Ann.\vol 304\pages 21--51\yr1996\endref

\ref\key CS \by C. De Concini and E. Strickland\paper On the variety
of complexes \pages 57--77\yr1981\vol 41\jour Adv. in Math.\endref

\ref\key EGA-I\by A. Grothendieck and J. Dieudonn{\'e}\book
{\'E}l{\'e}ments de G{\'e}om{\'e}trie
alg{\'e}brique I\bookinfo Grundlehren 166\publaddr
Berlin-Heidelberg-New York: Springer\yr1971\endref

\ref \key Fa \by G. Faltings\paper Moduli-stacks for bundles on semistable
curves\pages 489--515\yr1996\vol 304\jour Math. Ann.\endref
 
\ref\key FH \by W. Fulton and J. Harris\book Representation Theory : 
A first course\yr1991\bookinfo Graduate Texts in Mathematics
\vol 129\publaddr Springer-Verlag New York Inc.\endref

\ref\key Fu \by W. Fulton\paper Flags, Schubert polynomials, degeneracy 
loci and determinantal formula \yr1992\jour Duke Math. J.\vol 65\pages
381--420\endref

\ref \key Kn \by F. Knop\paper Der kanonische Moduleines Invariantenrings \pages
40--54\yr1989\vol 127\jour Joural of algebra\endref

\ref \key La \by H. Lange\paper Universal families of extensions \pages
101--112\yr1983\vol 83\jour Joural of algebra\endref

\ref\key Ne \by P.E.Newstead\book Introduction to
moduli problems
and orbit spaces\yr1978\bookinfo TIFR lecture
notes\publaddr New
Delhi: Narosa\endref

\ref\key NR \by M.S.Narasimhan and
T.R. Ramadas\paper Factorisation of 
generalised theta functions I\pages
565--623\yr1993\vol 114 \jour
Invent. Math.\endref

\ref\key NS \by D.S.
Nagaraj and C.S. Seshadri\paper Degenerations
of the moduli spaces of
vector bundles on curves I\pages 101--137
\yr1997\vol 107 \jour Proc.
Indian Acad. Sci.(Math. Sci.)\endref

\ref\key Pa \by C. Pauly\paper
Espaces modules de fibr{\'e}s
paraboliques et blocs conformes\pages
217--235\yr1996\vol 84 \jour
Duke Math.\endref

\ref\key Ra \by T.R.
Ramadas\paper Factorisation of 
generalised theta functions II: the
Verlinde formula\pages
641--654\yr1996\vol 35
\jour
Topology\endref

\ref\key Se \by C.S. Seshadri\paper Fibr{\'e}s
vectoriels sur les courbes
alg{\'e}briques\vol 96\jour Ast{\'e}risque\yr
1982\endref

\ref \key Si \by C. Simpson \paper Moduli of representations
of the
fundamental group of a smooth projective variety I
\pages
47--129\vol 79 \yr1994\jour I.H.E.S. Publications
Math{\'e}matiques\endref

\ref\key Su \by X. Sun\paper Degeneration of moduli spaces and 
generalized theta functions
\jour J. of Algebric geomety (to appear)\endref

\ref\key Sw \by R.G. Swan\paper On seminormality \pages 210--229
\yr1980\vol 67\jour J. of Algebra\endref

\ref\key TB\by M. Teixidor I Bigas\paper Moduli of vector bundles
on treelike curves \pages 341--348
\yr1991\vol 290\jour Math. Ann.\endref

\ref\key Tra\by C. Traverso\paper Seminormality and Picard group
\pages 585--595
\yr1970\vol 24\jour Ann. Scuola Norm. Sup. Pisa\endref

\ref \key Tri \by V. Trivedi\paper The seminormality property of
circular
complexes \pages 227--230\yr1991\vol 101\jour Proc. Indian Acad.
Sci.
(Math. Sci.) \endref 






















\endRefs
\enddocument